\documentclass[11pt,reqno]{amsart}
\usepackage[utf8]{inputenc}
\usepackage{geometry}
\usepackage{xcolor}
\usepackage{units}
\usepackage{enumitem}
\usepackage{amstext}
\usepackage{amsthm}
\usepackage{amssymb}
\usepackage{graphicx}
\PassOptionsToPackage{normalem}{ulem}
\usepackage{ulem}
\usepackage[unicode=true,
 bookmarks=false,
 breaklinks=false,pdfborder={0 0 1},backref=false,colorlinks=false]
 {hyperref}

\makeatletter

\newcommand{\lyxdot}{.}

\providecolor{lyxadded}{rgb}{0.047,0.78,0.44}
\providecolor{lyxdeleted}{rgb}{1,0,0}
\DeclareRobustCommand{\mklyxadded}[1]{\textcolor{lyxadded}\bgroup#1\egroup}
\DeclareRobustCommand{\mklyxdeleted}[1]{\textcolor{lyxdeleted}\bgroup\mklyxsout{#1}\egroup}
\DeclareRobustCommand{\mklyxsout}[1]{\ifx\\#1\else\sout{#1}\fi}
\numberwithin{equation}{section}
\numberwithin{figure}{section}
\usepackage{amscd}
\usepackage{latexsym}

\usepackage{amscd}
\usepackage[mathscr]{euscript}
\usepackage{mathrsfs}
\usepackage{xcolor}
\usepackage{color}

\usepackage{lscape}%to rotate pages
\usepackage{epstopdf}

\usepackage{subcaption}  % Added for using subfigures

\usepackage{array}
\usepackage{tabularx}
\usepackage{longtable}
\usepackage{ltxtable}

\usepackage{color}

\usepackage{tikz}
\usetikzlibrary{decorations,arrows}
\usetikzlibrary{decorations.pathmorphing}
\usepgflibrary{decorations.pathreplacing} 
\usepackage{pgf}
\usetikzlibrary{patterns}

\tikzset{
  schraffiert/.style={pattern=horizontal lines,pattern color=#1},
  schraffiert/.default=black
}

\tikzset{
    ultra thin/.style= {line width=0.1pt},
    very thin/.style=  {line width=0.2pt},
    thin/.style=       {line width=0.4pt},% thin is the default
    semithick/.style=  {line width=0.6pt},
    thick/.style=      {line width=0.8pt},
    very thick/.style= {line width=1.2pt},
    ultra thick/.style={line width=2.4pt}
}

 \usetikzlibrary{shapes}

\definecolor{hellgrau}{rgb}{0.93,0.93,0.93}
\definecolor{hellergrau}{rgb}{0.97,0.97,0.97}
\definecolor{hellgruen}{rgb}{0.6,1.35,0.5}
\definecolor{grau}{rgb}{0.93,0.93,0.93}
\definecolor{hellblau}{rgb}{0.8,0.8,2.0}
\definecolor{blau}{rgb}{0.3,0.5,2.0}
\definecolor{hellrot}{rgb}{2.0,0.6,0.6}
\definecolor{gruen}{rgb}{0.3,0.75,0.2}
\definecolor{rot}{rgb}{0.9,0.1,0.1}
\definecolor{orang}{rgb}{1.3,0.65,0}

\allowdisplaybreaks

\DeclareMathAlphabet{\mathpzc}{OT1}{pzc}{m}{it}

\newcommand{\supp}{{\operatorname{supp}}}			%support
\allowdisplaybreaks

\usepackage{xcolor}
\newcommand{\Hmm}[1]{\leavevmode{\marginpar{\tiny%
$\hbox to 0mm{\hspace*{-0.5mm}$\leftarrow$\hss}%
\vcenter{\vrule depth 0.1mm height 0.1mm width \the\marginparwidth}%
\hbox to 0mm{\hss$\rightarrow$\hspace*{-0.5mm}}$\\\relax\raggedright #1}}}

\definecolor{green}{RGB}{0, 180, 0}
\definecolor{cyan}{RGB}{0, 180, 180}
\definecolor{yellow}{RGB}{211,211,0}

\renewcommand{\subset}{\subseteq}
\makeatother

\begin{document}
%%%%%%%%% canonical symbols %%%%%%%%%%%%%%%%%%%%%%%%%%

\global\long\def\ui{\mathbf{\textrm{i}}}%

\global\long\def\ue{\mathbf{\textrm{e}}}%

\global\long\def\ud{\mathbf{\textrm{d}}}%

\global\long\def\rmi{\mathbf{\textrm{i}}}%

\global\long\def\rme{\mathbf{\textrm{e}}}%

\global\long\def\rmd{\mathbf{\textrm{d}}}%

\global\long\def\sgn{\mathrm{sign}}%

\global\long\def\tr#1{\mathop{tr}\left(#1\right) }%

\global\long\def\id{\mathbf{1}}%

\global\long\def\C{\mathbb{C}}%

\global\long\def\R{\mathbb{R}}%

\global\long\def\Q{\mathbb{Q}}%

\global\long\def\N{\mathbb{N}}%

\global\long\def\Z{\mathbb{Z}}%

\global\long\def\ZZ{\mathbb{Z}_{2}}%

\global\long\def\Id{\mathbf{\mathbb{I}}}%

\global\long\def\V{\mathcal{V}}%

\global\long\def\set#1#2{\left\{  #1\thinspace:\thinspace#2\right\}  }%

%%%%%%%%% Paper specific symbols %%%%%%%%%%%%%%%%%%%%%%%%%%

\global\long\def\SU{\mathbf{U}}%

\global\long\def\p{\mathtt{p}}%
 % periodic orbit

\global\long\def\q{\mathtt{q}}%
% another periodic orbit

\global\long\def\PO#1{\mathcal{PO}_{#1} }%

\global\long\def\graph{\Gamma}%

\global\long\def\edges{\mathcal{E}}%

\global\long\def\vtcs{\mathcal{V}}%

\global\long\def\bvec{\underline{b}}%

\global\long\def\btvec{\underline{\widetilde{b}}}%

\global\long\def\supp{\mathrm{supp}}%

\global\long\def\pop{\mathrm{popcnt}}%

\global\long\def\Eubonds{\mathrm{\mathcal{B}_{\textrm{Euler}}}}%

\global\long\def\Eucount{\mathrm{N_{\textrm{Euler}}}}%

\global\long\def\bonds{\mathrm{\mathcal{B}}}%

\global\long\def\bbonds#1{\mathcal{B}_{\bvec}( #1) }%

\global\long\def\mat#1{Y_{#1} }%

\global\long\def\adj{A}%

\title{Spectral statistics of preferred orientation quantum graphs}
\author{Ram Band$^{1}$, Pavel Exner$^{2,3}$, Divya Goel$^{4}$ and Aviya
Strauss$^{5}$}
\address{$^{1}$Department of Mathematics, Technion--Israel Institute of Technology,
Haifa, Israel}
\address{$^{2}$Doppler Institute, Czech Technical University, Prague, Czechia}
\address{$^{3}$Nuclear Physics Institute, Czech Academy of Sciences, Řež,
Czechia}
\address{$^{4}$Department of Mathematical Sciences, Indian Institute of Technology
(BHU), Varanasi, 221005, India}
\address{$^{5}$Faculty of Electrical and Computer Engineering, Technion--Israel
Institute of Technology, Haifa, Israel}
\begin{abstract}
We study the spectral statistics of quantum (metric) graphs whose
vertices are equipped with preferred orientation vertex conditions.
When comparing their spectral statistics to those predicted by suitable
random matrix theory ensembles, one encounters some deviations. We
point out these discrepancies and demonstrate that they occur in various
graphs and even for Neumann-Kirchhoff vertex conditions, which was
overlooked so far. Detailed explanations and computations are provided
for this phenomena. To achieve this, we explore the combinatorics
of periodic orbits, with a particular emphasis on counting Eulerian
cycles.
\end{abstract}

\maketitle

\section{Introduction\protect\label{sec:=000020Introduction}}

Since the seminal paper of Kottos and Smilansky \cite{KotSmi_prl97}
it is known that quantum graphs \cite{BerKuc_graphs,Berkolaiko_qg-intro17,KosNic_book22,BanGnu_qg-exerices18,GnuSmi_ap06,Kurasov_book24}
are a suitable class of systems on which chaotic properties can be
tested. One of the accepted conclusions is that the spectral statistics
of quantum graphs fall into several universal classes, in particular,
that the \emph{Gaussian Orthogonal} \emph{Ensemble} (GOE) distribution
can be observed only if the system exhibits invariance with respect
to the time reversal. One aim of this paper is to challenge this `rule'
by providing examples of graphs having the GOE eigenvalue statistics
despite being time-reversal asymmetric.

To give an example, one may consider a quantum particle living on
an \emph{octahedron} of incommensurate edge lengths assuming that
the boundary values of the wavefunctions and their derivatives at
each vertex are matched through the conditions,
\begin{equation}
\psi_{j+1}-\psi_{j}+i(\psi'_{j+1}+\psi'_{j})=0,\;\quad j=1,2,3,4\;(\mathrm{mod}\,4),\label{bc_comp}
\end{equation}
which are obviously non-invariant with respect to complex conjugation
that represents time reversal. 

We denote by $\left\{ k_{n}\right\} _{n=1}^{\infty}$ the square roots
of the eigenvalues of such a quantum graph (see a detailed description
of the model in Section~\ref{sec:=000020The-model}). First, we consider
the nearest-neighbour spacing distribution which is given by
\begin{equation}
P(x)=\lim_{N\rightarrow\infty}\frac{1}{N}\sum_{i=1}^{N}\delta(x-(k_{i+1}-k_{i})).\label{eq:=000020neares=000020neighbours=000020def}
\end{equation}

\begin{figure}[h!]
\centering \includegraphics[width=0.6\textwidth]{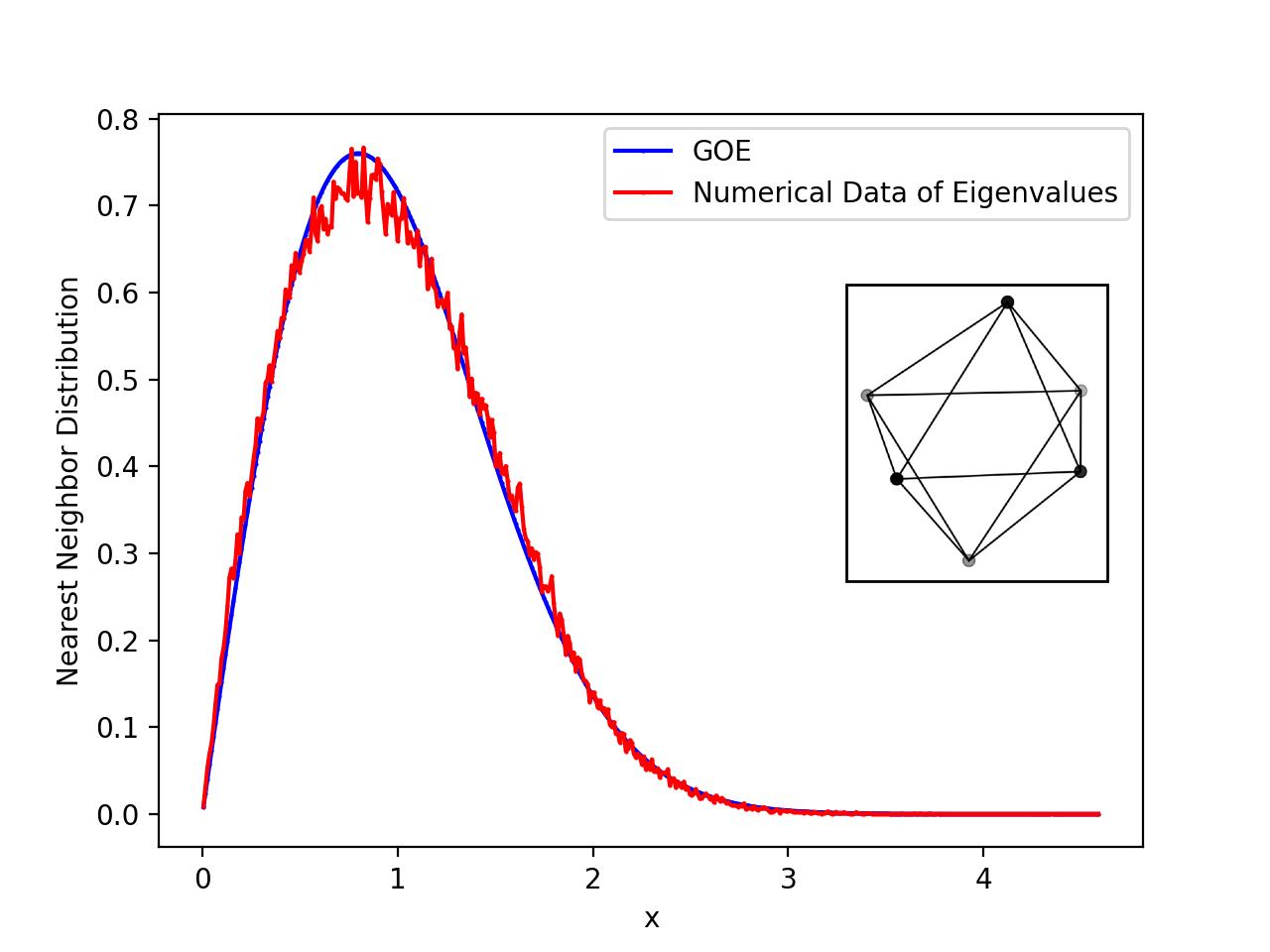}
\caption{The nearest-neighbour distribution of the first $4\cdot10^{5}$ eigenvalues
of an octahedron graph with incommensurate edge lengths and preferred
orientation vertex conditions, (\ref{bc_comp}). \protect\label{fig:=000020Nearest=000020neighbours=000020-=000020Octahedron}}
\end{figure}

The nearest-neighbour distribution of an octahedron graph with incommensurate
edge lengths follows the Gaussian Orthogonal Ensemble (GOE), as is
shown in Figure~\ref{fig:=000020Nearest=000020neighbours=000020-=000020Octahedron}.
% -------------- %
 The condition (\ref{bc_comp}) and its extensions have further interesting
consequences for the nearest-neighbour distributions, as we will show
in Section~\ref{sec:=000020Numerics=000020-=000020Nearest=000020Neighbours}
and explain in Section~\ref{sec:=000020Discussion=000020of=000020Nearest=000020neighbours}.

\begin{figure}[h!]
\centering \includegraphics[scale=0.54]{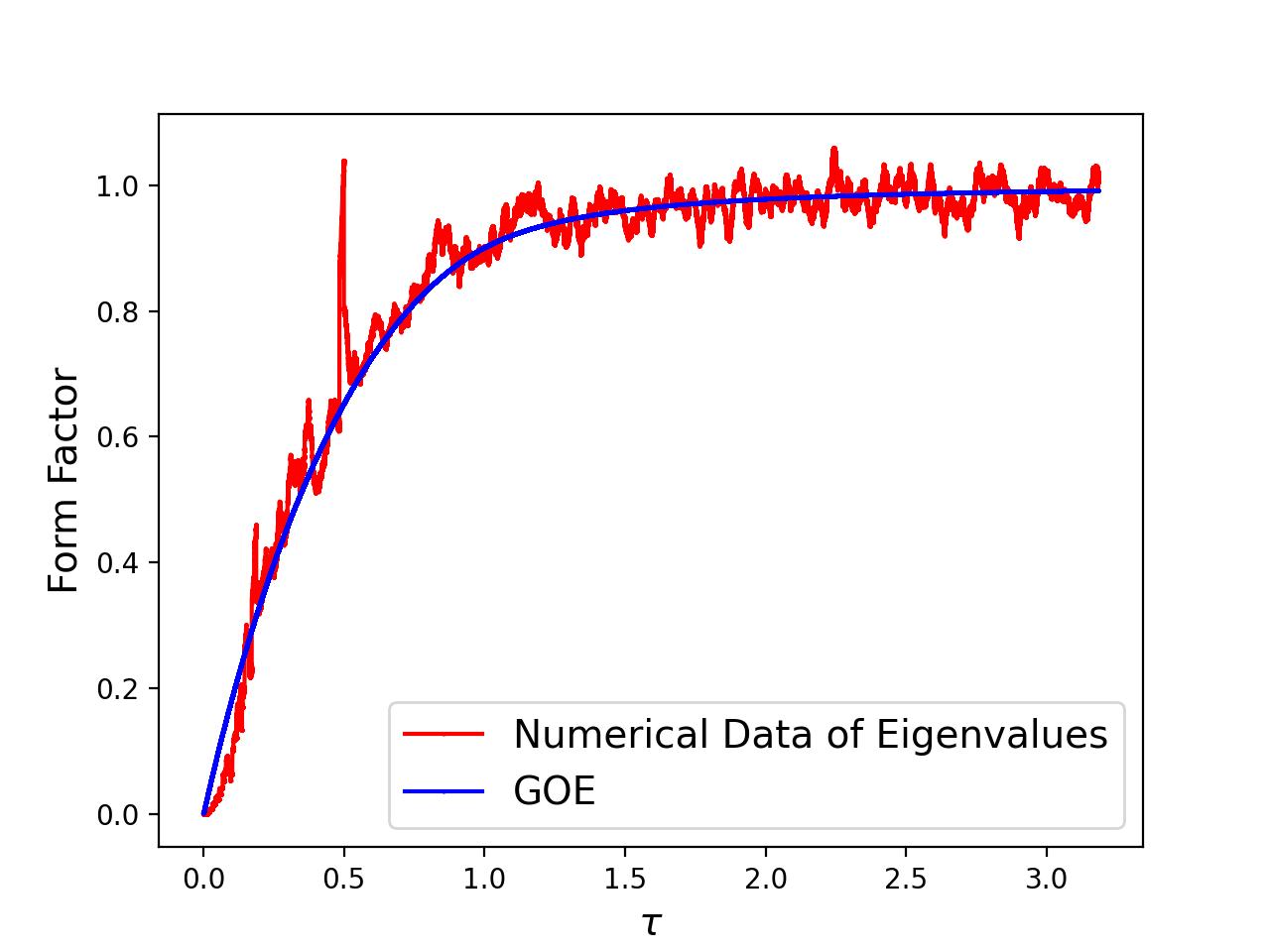} \caption{The form factor of the octahedron graph with incommensurate edge lengths
and preferred orientation vertex conditions, (\ref{bc_comp}). The
numerics was done using eigenvalues no. $5\cdot10^{4}-1.5\cdot10^{5}$.
\protect\label{fig:=000020Form=000020factor=000020of=000020Octahedron}}
\end{figure}

Furthermore, when going beyond nearest-neighbour distribution, one
considers the two-point correlation function, 
\begin{equation}
R_{2}(x)=\lim_{N\rightarrow\infty}\frac{1}{N}\sum_{i=1}^{N}\sum_{j=1}^{N}\delta(x-(k_{j}-k_{i})).\label{eq:=000020two-point=000020correlation}
\end{equation}
The Fourier transform of $R_{2}(x)$ is called the form factor, 
\begin{equation}
K(\tau)=\int_{-\infty}^{\infty}\rme^{2\pi\rmi x\tau}(R_{2}(x)-1)\thinspace\rmd x.\label{eq:=000020form=000020factor}
\end{equation}
From a global viewpoint, the form factor of the same octahedron graph
also exhibits GOE like behaviour. Nevertheless, there is a clear deviation
from GOE in the form of a sharp peak at $\tau=1/2$, as is shown in
Figure~\ref{fig:=000020Form=000020factor=000020of=000020Octahedron}.

To the best of our knowledge, this specific deviation from the form
factor was never mentioned before in the literature. We explain this
phenomenon in Section~\ref{sec:=000020form=000020factor}, by providing
the required analysis of the form factor. Additionally, we discuss
the extensions of these phenomena to other graphs and other vertex
conditions. 

In the next section, we describe in more detail the model and the
background in random matrix theory.

\section{The model and some background\protect\label{sec:=000020The-model}}

\subsection{Graphs with preferred orientation conditions}

The motion on the graph edges is free, being described by the Laplacian,
$\psi_{j}\mapsto-\psi_{j}''$; the nontrivial part comes from the
conditions matching the wave functions at the vertices. We write the
boundary values at each vertex $v$ as columns, $\Psi(v):=\{\psi_{j}(v)\}$
and $\Psi'(v):=\{\psi'_{j}(v)\}$, understood as limits at the endpoint;
then the most general way to make the Laplacian a self-adjoint operator
is to require % ------------- %
 
\begin{equation}
(U-I)\Psi(v)+\rmi(U+I)\Psi'(v)=0,\label{bc_matrix}
\end{equation}
% ------------- %
where $U$ is an auxiliary unitary $d\times d$ matrix ($d$ being
the degree of the vertex). The origin of this formulation of the vertex
conditions is usually referred to \cite{KosSch_jpa99} but in fact
they appeared already in \cite{RofeBeketov1969}.

We work with a simple preferred-orientation coupling proposed in \cite{Exner2018a}
in which 
\begin{equation}
U=\left(\begin{array}{cccccc}
0 & 1 & 0 & \ldots & \ldots & 0\\
0 & 0 & 1 & \ddots &  & \vdots\\
\vdots &  & 0 & \ddots & 0 & \vdots\\
\vdots &  &  & \ddots & 1 & 0\\
0 & \ldots &  &  & 0 & 1\\
1 & 0 & \ldots &  & \ldots & 0
\end{array}\right)\label{eq:=000020circulant=000020matrix}
\end{equation}
is the \emph{circulant matrix} \cite{Davis1979}. As an example, writing
\eqref{bc_matrix} for a degree four vertex in components we get the
conditions \eqref{bc_comp}. The on-shell\footnote{Where if one wants to choose a different length scale (i.e., change
units), it is enough to replace $k$ by $k\ell$ for a fixed value
$\ell>0$.} scattering matrix at a vertex described by the conditions \eqref{bc_matrix}
is % ------------- %
\begin{equation}
S(k)=\frac{k-1+(k+1)U}{k+1+(k-1)U}\,;\label{eq:=000020S(k)}
\end{equation}
% ------------- %
 It is easy to check that the S-matrix \eqref{eq:=000020S(k)} is
\emph{not invariant with respect to transposal}, which implies in
our case that the transport through the vertex \emph{not} being time-reversal
invariant \cite{Exner2021}. 

To complement the description above, we mention the case of Neumann-Kirchhoff
vertex conditions, in which the scattering matrix $S(k)$ is independent
of $k$ and its entries equal $\left[S(k)\right]_{i,j}=2/d-\delta_{i,j}$
for a vertex of degree $d$. Later in the paper we compare results
obtained with preferred orientation vertex conditions with results
obtained for Neumann-Kirchhoff vertex conditions.

\subsection{Spectral statistics of quantum graphs and random matrix theory -
existing results}

A driving force in connecting the spectral statistics of chaotic systems
to random matrix theory (RMT) lies in the Bohigas-Giannoni-Schmit
(BGS) conjecture \cite{BohGiaSch_prl84}. In the realm of quantum
graphs the conjecture says that the spectral statistics of graphs
with incommensurate edge lengths and ``sufficient connectivity'' exhibit
universal statistical properties governed by random matrix theory.
Moreover, the symmetry class to which the graphs belong dictates which
RMT ensemble describes their spectral statistics. This is still somewhat
vague phrasing and substantial effort was made in order to even formalize
this conjecture and specify the exact conditions for its validity
(both for quantum graphs and for other system within quantum chaos).
The connection between quantum graph spectral statistics and RMT was
first investigated by Kottos and Smilansky \cite{KotSmi_prl97,KotSmi_ap99},
and later extended in works with Schanz \cite{SchSmi_pmb00,SchSmi_ejc01,KotSch_phys01},
where spectral fluctuations were analyzed using combinatorial methods
(see also \cite{GnuSmi_ap06} for an extensive review of these results
and the ones which follow). Barra and Gaspard \cite{BarGas_jsp00}
contributed a careful study of the nearest-neighbour spacing distribution,
helping to establish statistical links with RMT. Tanner \cite{Tan_jpa00,Tan_jpa01}
introduced unitary-stochastic matrix ensembles for classifying when
graphs display RMT-like behavior. In particular, he conjectured that
the RMT-like behavior is observed if the spectral gap of a certain
transition matrix closes slow enough.Berkolaiko and Keating explored
spectral statistics in star graphs (where the aforementioned condition
is violated) \cite{BerKea_jpa99,Berkolaiko_thesis00}, and Berkolaiko,
in part with Schanz and Whitney, refined these methods analyzing form
factors via periodic orbit expansions and diagrammatic approaches
\cite{BerSchWhi_prl02,BerSchWhi_jpa03,Ber_wrm04,Ber_incol06}. Bolte
and Harrison studied spectral statistics for the spin-orbit coupling
and for the Dirac operator on graphs \cite{BolHar_jpa03,BolHar_jpa03a,BolHar_incol06}.

Additional analysis have been obtained using field-theoretic and supersymmetric
approaches. In particular, Gnutzmann and Altland \cite{GnuAlt_prl04,GnuAlt_pre05}
applied the nonlinear sigma model to show that spectral correlations
of individual quantum graphs match RMT predictions. These techniques
were further developed by Pluha\v{r} and Weidenm\"uller \cite{PluWei_prl13,PluWei_prl14,PluWei_jpa15,Weidenmueller_chapter20},
who established RMT universality using diagrammatic and supersymmetric
formulations.

Further theoretical works have continued refining the boundary between
the classical connection to RMT symmetry classes and its breakdown.
Joyner, M\"uller, Sieber \cite{JoyMueSie_epl14}, as well as Akila
and Gutkin \cite{AkiGut_jpa15,Akiut_jpa19}, showed that systems without
spin can still yield GSE-type spectral statistics. Harrison, Swindle
and Winn showed intermediate spectral statistics for various models
\cite{HarSwi_jpa19,HarWin_jpa12}. Band, Harrison, Hudgins, Joyner
and Sepanski studied the the variance of coefficients of the characteristic
polynomial of the quantum evolution operator \cite{BanHarJoy_jpa12,BanHarSep_jmaa19,HarHud_epl22,HarHud_jpa22}
(appearing already in the earlier works of Tanner mentioned above).
Most recently, Gnutzmann and Smilansky \cite{GnuSmi_jpa24} emphasized
that RMT-like spectral statistics do not necessarily indicate chaotic
classical dynamics.

On the experimental side, realizations of quantum graphs have played
a crucial role both in validating RMT connections and in showing deviations
from them \cite{OleSir_pre11,FarAkhLawBiaSir_pre24,CheGluKohGuhDie_arXiv25,DieKlaMasMisRicSkiWun_pre24,HulBauPakSavZycSir_pre04,DieYunBiaBauLawSir_pre17,HoLuKuSto_pre21,LuHofKuhSto_ent23,RehAllJoyMueSieKuhSto_prl16}.

\section{Numerical results for nearest neighbour distribution\protect\label{sec:=000020Numerics=000020-=000020Nearest=000020Neighbours}}

Together with the numerical results shown in the introduction (Figure~\ref{fig:=000020Nearest=000020neighbours=000020-=000020Octahedron})
for the octahedron graph we also consider other graphs and/or couplings
in order to better elucidate the mechanism responsible for the observed
effects. First of all, let us note that the eigenvalue counting function
of a finite graph with arbitrary self-adjoint vertex conditions satisfies
Weyl's law, $N(k)=\frac{L}{\pi}k+\mathcal{O}(1)$ as $k\to\infty$,
where $L=\sum_{j}\ell_{j}$ is the sum of all the edge lengths. This
is known theoretically first from \cite{KotSmi_ap99} and later in
\cite[prop.~4.2]{BolEnd_ahp09} for the most general self-adjoint
conditions (including the preferred orientation conditions which we
consider here). As a consequence, the unfolding is trivial; the proper
scale to display the eigenvalue spacing is given by the simply scaled
momentum variable $\frac{L}{\pi}k$.

\begin{figure}[h!]
\centering \includegraphics[width=0.5\textwidth]{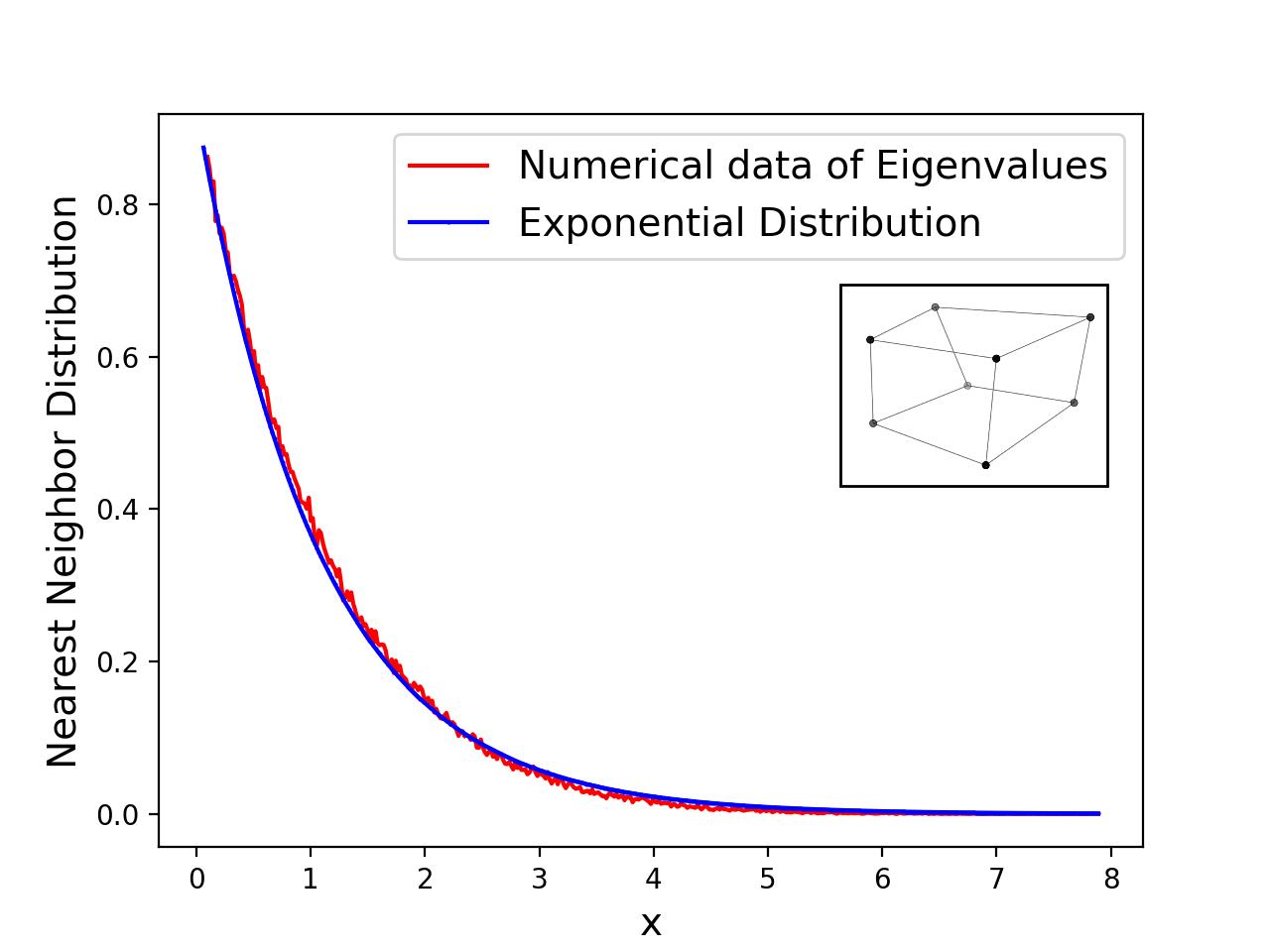}
\caption{The nearest-neighbour distribution of the first $2\cdot10^{5}$ eigenvalues
of a cube graph with incommensurate edge lengths and preferred orientation
vertex conditions.\protect\label{fig:=000020cubePoisson}}
\end{figure}

We return to Figure~\ref{fig:=000020Nearest=000020neighbours=000020-=000020Octahedron}
which shows that the nearest neighbour distribution of an octahedron
graph is GOE. Modifying the vertex degree can change the spectral
picture completely; to illustrate that, we show in Figure~\ref{fig:=000020cubePoisson}
the eigenvalue spacing distribution for a cube graph with incommensurate
edge lengths and the same vertex coupling \eqref{bc_comp}. In this
case we observe that the distribution is of the Poisson type.

% -------------- %

\begin{figure}[h!]
\centering \includegraphics[width=0.5\textwidth]{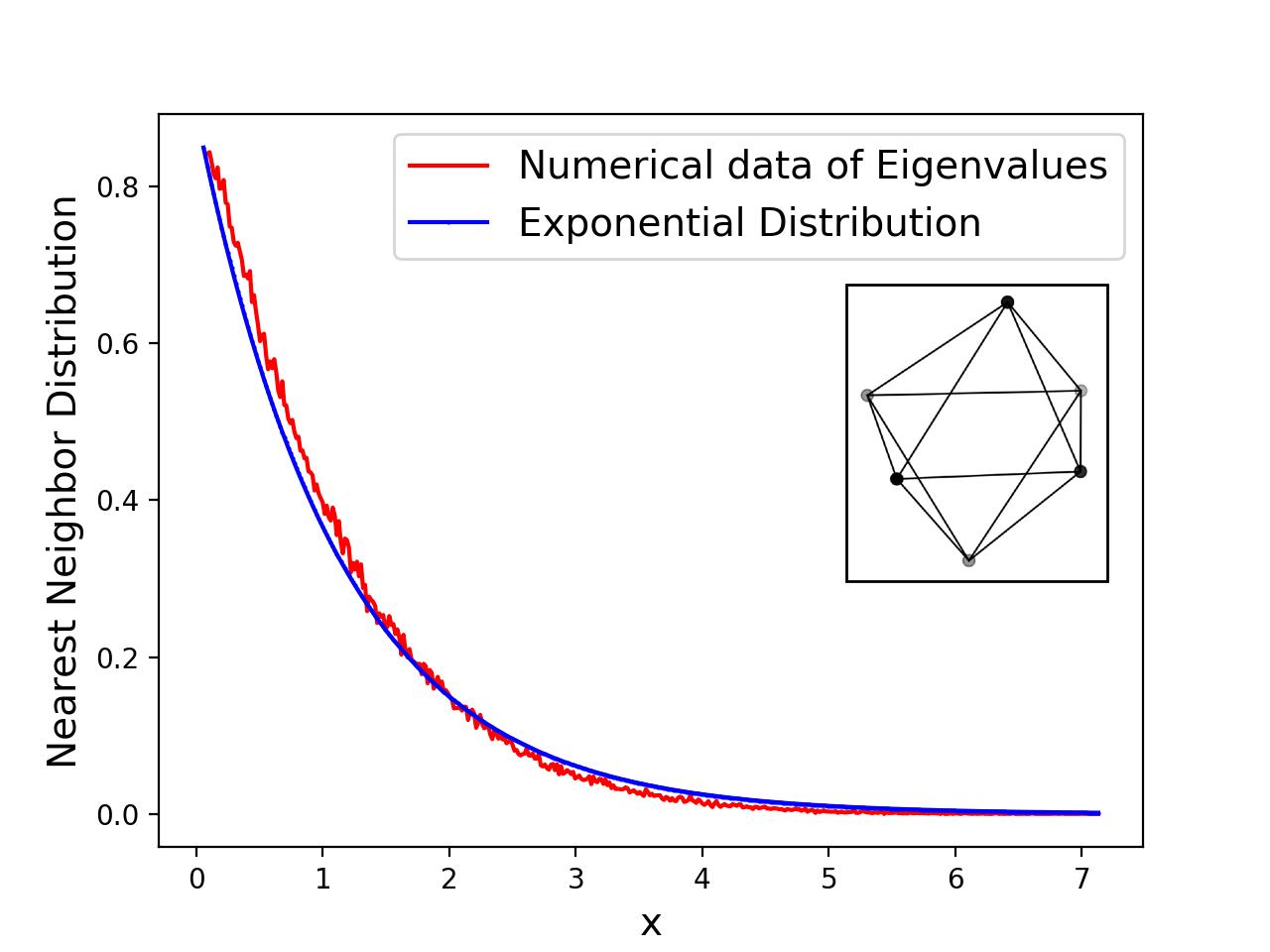}
\caption{The nearest-neighbour distribution of the first $2\cdot10^{5}$ eigenvalues
of an octahedron graph with incommensurate edge lengths and distorted
preferred orientation vertex conditions, (\ref{eq:=000020distorted_R})
for $\mu=0.01$.}
\label{fig:=000020distort_octahedron}
\end{figure}

Furthermore, it is not so much the vertex degree that determines the
statistics type, but rather the presence or absence of the eigenvalue
$-1$ in the spectrum of the matrix $U$ in (\ref{eq:=000020circulant=000020matrix}).
To illustrate this claim, consider the octahedron again, but replace
now the coupling \eqref{eq:=000020circulant=000020matrix} with a
`distorted' one referring to the modified matrix
\begin{equation}
U=\mathrm{e}^{\rmi\mu}\left(\begin{array}{cccc}
0 & 1 & 0 & 0\\
0 & 0 & 1 & 0\\
0 & 0 & 0 & 1\\
1 & 0 & 0 & 0
\end{array}\right),\label{eq:=000020distorted_R}
\end{equation}
for some $\mu>0$. Figure~\ref{fig:=000020distort_octahedron} shows
that even if the parameter $\mu$ is small, the modification changes
the picture completely; instead of the GOE we have the Poisson distribution\footnote{To be exact, in the numerics the distribution becomes Poissonian if
we consider a sufficiently wide energy interval; before we reach this
regime one observes a mixture between GOE and Poisson. }.

\section{Discussion of the nearest neighbour distribution\protect\label{sec:=000020Discussion=000020of=000020Nearest=000020neighbours}}

\subsection{High energy asymptotic of the unitary evolution operator\protect\label{subsec:=000020High-energy-asymptotic}}

The occurrence of Poisson distribution in Figures~\ref{fig:=000020cubePoisson}
and \ref{fig:=000020distort_octahedron} is not surprising. In both
cases the eigenvalue $-1$ is missing in the spectrum of the matrix
$U$ determining the coupling, and as a result, we get $\lim_{k\to\infty}S(k)=I$
from \eqref{eq:=000020S(k)}. This means that the considered graph
turns at high energies effectively into a union of disconneted edges
with Neumann endpoints.

The existence of the GOE distribution in Figure~\ref{fig:=000020Nearest=000020neighbours=000020-=000020Octahedron}
might be less obvious but to understand it one has to realize that
while the time-reversal symmetry is violated for any $k$, the degree
of the violation varies. We refer to \cite[eq.~(5)]{Exner2018a},
in whichthe entries of the matrix (\ref{eq:=000020S(k)}) were expressed
as
\begin{equation}
S_{i,j}(k)=\frac{1-\eta^{2}}{1-\eta^{d}}\left\{ -\eta\frac{1-\eta^{d-2}}{1-\eta^{2}}\delta_{i,j}+(1-\delta_{i,j})\eta^{(j-i-1)\mod d}\right\} ,\label{eq:=000020S(k)=000020using=000020eta}
\end{equation}
where $\eta:=\frac{1-k}{1+k}$. From here, one can check that the
S-matrix of a vertex of even degree $d$ has in the high-energy limit
the entries % -------------- %
\begin{equation}
S_{i,j}=-\frac{2}{d}(-1)^{i+j}+\delta_{i,j},\label{eq:=000020Even=000020degree=000020asymptotics}
\end{equation}
which differ from their Neumann-Kircchoff counterparts only by sign
on the diagonal and the ``even'' (i.e., even $i+j$) off-diagonals. 

If the vertex is of odd degree then at high energy asymptotics there
is an effective decoupling of this vertex into $d$ disjoint vertices
of degree one each with Neumann-Kirchhoff conditions, i.e., $S_{i,j}=\delta_{i,j}$.

In both cases (even and odd degrees) the limiting matrix is transpose
invariant. Hence, even though the time reversal invariance is violated
at any finite energy, it does survive asymptotically. For an even
$d$, in addition, the edges remain coupled asymptotically, which
is the reason for the GOE statistics for graphs with even degree vertices
and preferred orientation vertex conditions (such as the octahedron
considered in Figure~\ref{fig:=000020Nearest=000020neighbours=000020-=000020Octahedron}.
We supplement this observation with a quantitative discussion of the
time reversal invariance.

\subsection{On measuring the time reversal invariance violation}

The natural measure of time-reversal invariance violation says how
much the S-matrix differs from its transpose. We are able to express
the violation measure quantitatively by defining %----------------%
\begin{equation}
\mathcal{M}(k):=\|S(k)-S(k)^{T}\|\,;\label{eq:=000020measure-1}
\end{equation}
from the unitarity of $S(k)$ and the triangle inequality we conclude
that the norm cannot exceed two. We employ \eqref{eq:=000020S(k)=000020using=000020eta}
and write the matrix elements of $M(k)=S(k)-S(k)^{T}$ at a vertex
of degree $d$ as
\begin{equation}
M_{ij}(k)=\frac{1-\eta^{2}}{1-\eta^{d}}\,\big\{\eta^{(j-i-1)\mod d}-\eta^{(i-j-1)\mod d}\big\}.\label{eq:=000020M=000020matrix}
\end{equation}
To indicate the dependence on the vertex degree, we add the index
$d$ to the symbol $\mathcal{M}$, in what follows. For convenience
we will consider separately the prefactor $P_{d}=\frac{1-\eta^{2}}{1-\eta^{d}}$
and the matrix $T$, whose entries are given by the curly bracket.
Being a difference of circulant matrices, $T$ is also circulant which
means that its eigenvalues can be expressed explicitly,
\begin{equation}
\lambda_{n}(\eta)=\sum_{j=1}^{d-1}\big(\eta^{j-1}-\eta^{d-j-1}\big)\,\omega^{jn},\quad n=0,\dots,d-1\,,\label{eq:=000020eigenvalues}
\end{equation}
where $\omega:=\mathrm{e}^{2\pi i/d}$, and in particular, we have
$\lambda_{0}=0$. The norm of $T$ is naturally $\|T\|=\max_{0\le n\le d-1}|\lambda_{n}(\eta)|$.

Being primarily interested in the high-energy behavior of $\mathcal{M}_{d}(k)$,
we have to distinguish the even and odd values of $d$. The prefactor
equals
\begin{equation}
P_{2m}(\eta)=\frac{1}{1+\eta^{2}+\cdots+\eta^{2m-2}}\,,\;\;P_{2m+1}(\eta)=\frac{1+\eta}{1+\eta+\cdots+\eta^{2m}},\label{eq:=000020prefactor}
\end{equation}
which, translated to the original momentum variable,implies 
\begin{equation}
P_{2m}(k)=\frac{1}{m}+\mathcal{O}(k^{-1})\quad\text{and}\quad P_{2m+1}(\eta)=\frac{2}{k}+\mathcal{O}(k^{-2}),\label{eq:=000020prefasympt}
\end{equation}
as $k\to\infty$. On the other hand, putting $\eta=-1+\delta$ we
get from \eqref{eq:=000020eigenvalues} the following asymptotic expansion,
\begin{equation}
\lambda_{n}(\eta)=\begin{cases}
2i\delta\sum_{j=1}^{(d-2)/2}(-1)^{j+1}(d-2j)\sin\frac{2\pi nj}{d}+\mathcal{O}(\delta^{2}) & d\;\text{even},\\
4i\sum_{j=1}^{(d-1)/2}(-1)^{j+1}\sin\frac{2\pi nj}{d}+\mathcal{O}(\delta) & d\;\text{odd},
\end{cases}\label{eq:=000020evasympt}
\end{equation}
where we used the fact that the members of the sum \eqref{eq:=000020eigenvalues}
appear in pairs having the same real parts while the imaginary ones
differ by sign. Using next the fact that $\eta=-1+2k^{-1}+\mathcal{O}(k^{-2})$,
we conclude from \eqref{eq:=000020measure-1}, \eqref{eq:=000020prefasympt}
and \eqref{eq:=000020evasympt} that $\mathcal{M}_{d}(k)=\mathcal{O}(k^{-1})$
holds as $k\to\infty$ for any natural $d\ge3$. Note that the asymptotic
behavior of both the prefactor and the eigenvalues depends on the
vertex parity, however, the differences compensate mutually in the
result. In the same way one can check that the violation measure $\mathcal{M}_{d}(k)=\mathcal{O}(k)$
as $k\to0$.

For low values of $d$ it is easy to evaluate the non-invariance measure
explicitly. In particular, the eigenvalues \eqref{eq:=000020eigenvalues}
are $\{0,\pm\sqrt{3}i(1-\eta)\}$ and $\{0,0,\pm2i(1-\eta^{2})\}$
for $d=3,4$, respectively, which yields 
\begin{equation}
\mathcal{M}_{3}(k)=\frac{4k\sqrt{3}}{3+k^{2}}\quad\text{and}\quad\mathcal{M}_{4}(k)=\frac{4k}{1+k^{2}}.\label{eq:=000020low_n}
\end{equation}
As a concluding remark we note that \eqref{eq:=000020low_n} shows
that $\mathcal{M}_{3}(k)$ saturates the unitarity bound $\mathcal{M}(k)\le2$
at $k=\sqrt{3}$, while $\mathcal{M}_{4}(k)$ does the same at $k=1$.

\section{The form factor\protect\label{sec:=000020form=000020factor}}

We focus in this section on the form factor, as given in (\ref{eq:=000020form=000020factor}).
We change the point of view from the the spectral statistics of the
graph eigenvalues to those of the eigenphases of the corresponding
unitary operator. This is a common practice within the spectral theory
of quantum graphs (see \cite{BerWin_tams10} for justifications and
proofs). Indeed, one observes that the eigenphase form factor follows
the one of the eigenvalues (as demonstrated in Figure~\ref{fig:=000020Form=000020factor=000020-=000020complete=000020graphs},
as well as in the other figures in this section). To be specific,
denoting by $E$ the number of graph edges, the form factor of the
unitary operator is defined at the discrete times $\tau\in\frac{1}{2E}\N$
by (see e.g. \cite{Ber_incol06}) 
\begin{equation}
K_{\SU}(\tau)=\frac{1}{2E}\lim_{\Lambda\rightarrow\infty}\frac{1}{2\Lambda}\int_{-\Lambda}^{\Lambda}\left|\tr{\SU(k)}^{2E\tau}\right|^{2}\rmd k,\label{eq:=000020form=000020factor=000020eigenphases}
\end{equation}
where $\SU(k)$ is a unitary $2E\times2E$ matrix (often called the
unitary evolution operator) given by 
\[
\SU(k)=\exp\left(\rmi kL\right)\mathbf{S}(k),
\]
with $L$ being a diagonal matrix which stores the (directed) edge
lengths and $\mathbf{S}(k)$ is the global scattering matrix of the
graph. The matrix $\mathbf{S}(k)$ is comprised from the local vertex
scattering matrices $S(k)$ (such as in \eqref{eq:=000020S(k)}).
Note that we use bold font (as in $\mathbf{S}(k)$ and $\SU(k)$)
for matrices of dimensions $2E\times2E$ defined on the whole graph,
to be distinguished from the local vertex matrices $S$, $U$ which
are $d\times d$, where $d$ is the degree of the vertex. Since the
spectral statistics are dominated by the high energy asymptotics (see
discussion in Section~\ref{sec:=000020Discussion=000020of=000020Nearest=000020neighbours})
we may replace the energy dependent scattering matrix $\mathbf{S}(k)$
with its high energy limit, which is comprised by the local vertex
scattering matrices given in \eqref{eq:=000020Even=000020degree=000020asymptotics}
(for an even degree vertex). Doing so, one sees that the only energy
($k$) dependence in \eqref{eq:=000020form=000020factor=000020eigenphases}
enters via the matrix $\exp\left(\rmi kL\right)$. Expanding $\tr{\SU(k)}^{2E\tau}$
as a sum of products and performing the integral give the following
useful expansion (see \cite[eq. (18)]{Ber_incol06})
\begin{equation}
K_{\SU}(\tau)=\frac{1}{2E}\sum_{\p,\q}\frac{(2E\tau)^{2}}{r_{\p}r_{\q}}A_{\p}A_{\q}^{*}\delta_{L_{\p},L_{\q}},\label{eq:=000020form=000020factor=000020-=000020orbit=000020expansion}
\end{equation}
where as above $2E\tau\in\N$. The sum in (\ref{eq:=000020form=000020factor=000020-=000020orbit=000020expansion})
is over pairs of periodic orbits consisting of $2E\tau$ edges and
denoted by $\p,\q$. The total metric length of an orbit $\p$ is
denoted by $L_{\p}$, its overall scattering amplitude (which is a
product of $\mathbf{S}$ entries corresponding to the orbit) is denoted
by $A_{\p}$. Furthermore, sometimes an orbit might be written as
a repetition of a shorter orbit; in such a case the repetition number
is denoted by $r_{\p}$ (if an orbit cannot be written as such repetition,
then $r_{\p}=1$). The expansion (\ref{eq:=000020form=000020factor=000020-=000020orbit=000020expansion})
is the starting point of the derivations in this section.

\subsection{The form factor in the limit of infinite complete graphs}

\begin{figure}[!h]
\raggedright
\begin{minipage}[t]{0.45\columnwidth}%
\includegraphics[width=1.12\textwidth]{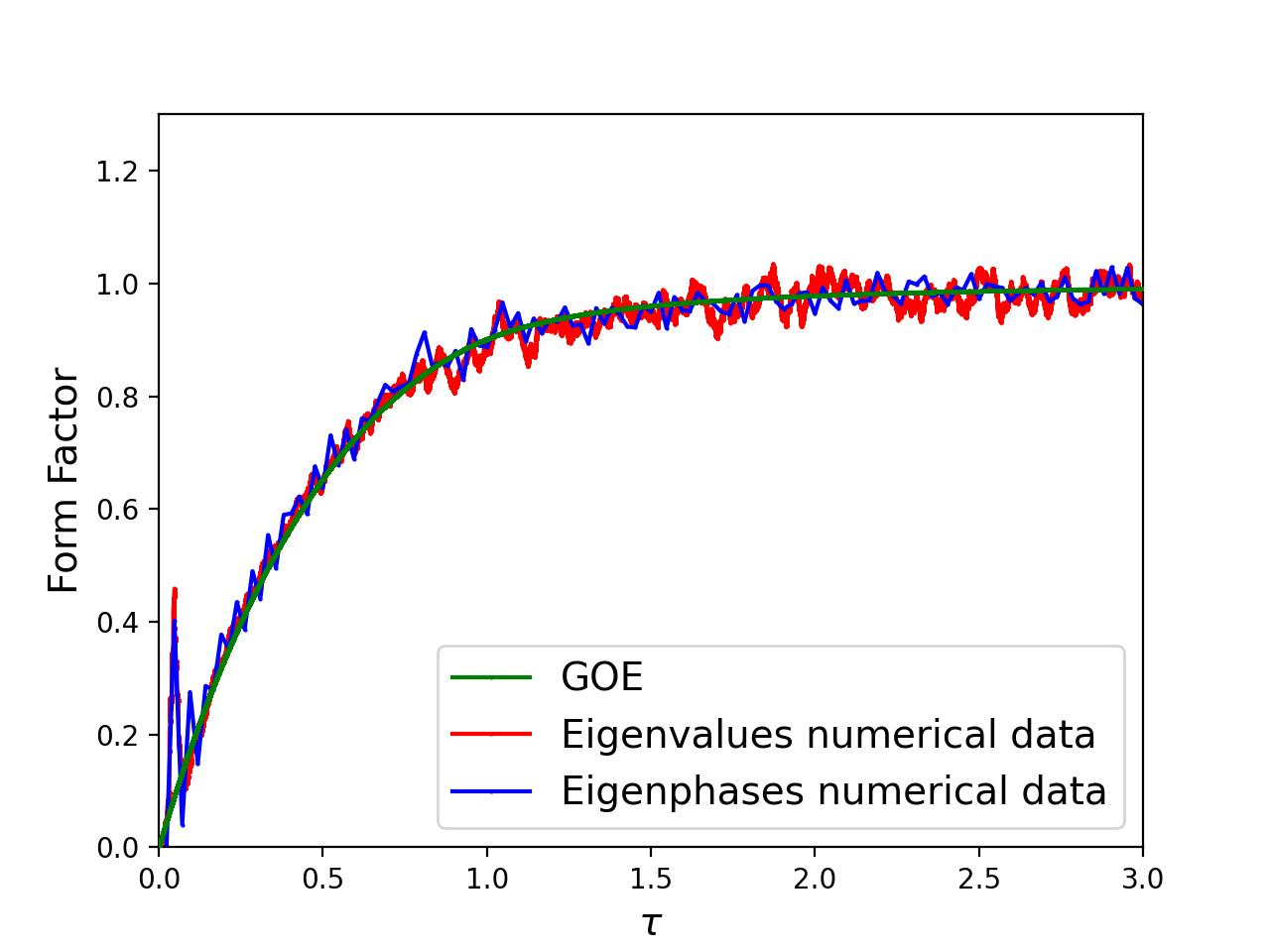}%
\end{minipage}$\quad\quad\quad$%
\begin{minipage}[t]{0.45\columnwidth}%
\includegraphics[width=1.12\textwidth]{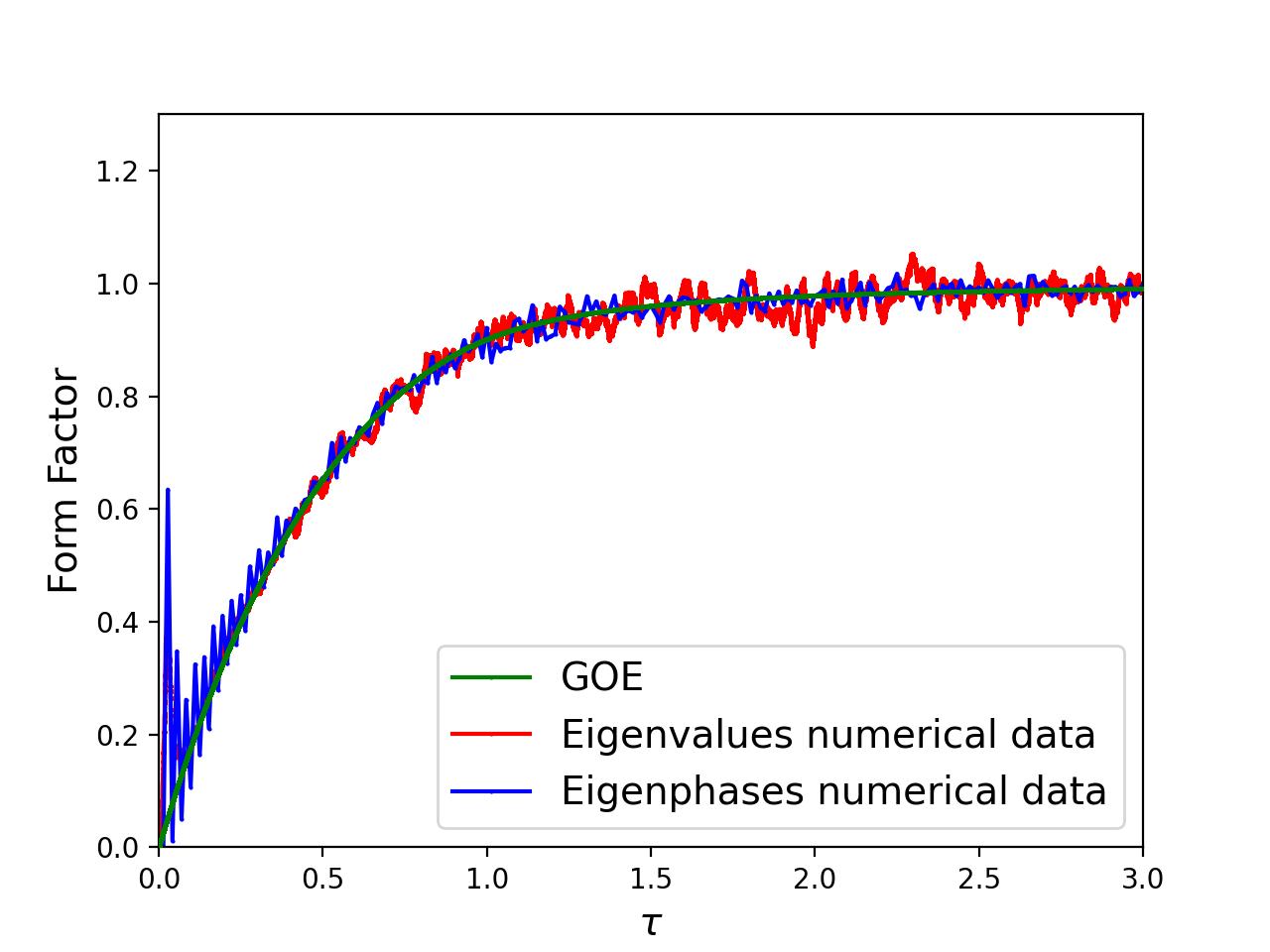}%
\end{minipage}\caption{The form factor of complete graphs (Left: $K_{7}$. Right: $K_{9}$).
The curves correspond to: numerics of eigenvalues (red), numerics
of the eigenphases (blue) and the theoretical GOE (green). \protect\label{fig:=000020Form=000020factor=000020-=000020complete=000020graphs}}
\end{figure}

We consider here the infinite family of graphs, $\left\{ K_{V}\right\} _{V\in2\N+1}$,
i.e., the family of complete graphs with an odd number of vertices.
In general, in order to analytically prove any kind of RMT-like behaviour
for graphs, one should consider families of increasing graphs (for
finite graphs deviations might occur, as we indeed observe in this
work). We choose here the particular family $\left\{ K_{V}\right\} _{V\in2\N+1}$
thanks to the somewhat easier book-keeping of the periodic orbits
of complete graphs and since all of its vertices have even degrees
(following the discussion in Section~\ref{subsec:=000020High-energy-asymptotic}).
In the current subsection we calculate the leading term of the form
factor in the limit of increasing graphs of this family; we verify
that it is indeed the leading term of the GOE form factor (see (\ref{eq:=000020Form=000020Factor=000020-=000020Theoretic})).
As a demonstration we provide in Figure~\ref{fig:=000020Form=000020factor=000020-=000020complete=000020graphs}
the numerical calculation of the form factor for the graphs $K_{7}$
and $K_{9}$, two particular members of this graph family.

Consider a complete graph with $V$ vertices and $E={V \choose 2}$
edges. As mentioned above, we take odd values of $V$ and equip all
vertices with preferred orientation conditions. Therefore in the high
energy limit we get the vertex conditions as in (\ref{eq:=000020Even=000020degree=000020asymptotics}),
which we use in the following computation. 

First, we classify the periodic orbits which are of a specific size
$n$ (namely, periodic orbits which consist of $n$ edges). Each such
orbit $\p$ has some $0\leq t\leq n$ transmission scattering events
and $n-t$ reflection scattering events; its overall scattering amplitude
is therefore $A_{\p}=\pm\left(\frac{2}{V-1}\right)^{t}\left(1-\frac{2}{V-1}\right)^{n-t}$,
since the degree of each vertex is $V-1$. The ambiguity in the $\pm$
sign is due to the term $(-1)^{i+j}$ in (\ref{eq:=000020Even=000020degree=000020asymptotics}),
when $i\neq j$ (i.e., transmission scattering event) and we show
in the following that it does not affect the final result. 

We claim that the number of such orbits (with $t$ transmissions and
$n-t$ reflections) is $\frac{2}{n}{n \choose t}{V \choose 2}(V-2)^{t-2}$,
assuming that $t\geq2$; this combinatorial result is explained in
the following. First, one picks the ``starting'' directed edge of
the orbit, and there are $2E=2{V \choose 2}$ options for such a choice.
Then, one chooses which of the $n$ scattering events are the $t$
transmissions, which gives the binomial factor ${n \choose t}$. When
forming such an orbit, we know for each reflection event what is the
next directed edge (it is just the last edge with a reversed direction);
but for (almost) each transmission event we have $V-2$ possibilities
to choose the next edge (it can be any edge emanating from the current
vertex, apart from retracing along the last edge). Nevertheless, in
this process, the two transmission events which appear after all other
transmissions (but before some possible reflections) are uniquely
determined in a way which ensures that the orbit returns back to the
right ``starting'' edge. This gives a factor of $(V-2)^{t-2}$,
noting that we assumed $t\geq2$ (the case $t<2$ is explained below).
Then, we need to use a factor of $\frac{1}{n}$ since we consider
periodic orbits up to cyclic shifts (and an orbits consists of $n$
edges). Hence, multiplying all of the above we get a total number
of $\frac{2}{n}{n \choose t}{V \choose 2}(V-2)^{t-2}$ orbits which
consist of $t$ transmissions and $n-t$ reflections (counted up to
cyclic shifts). To complement this computation we check the case $t<2$.
First, it is easy to see that there are no periodic orbits (on the
complete graph) with only a single transmission ($t=1$). Hence, we
are left to deal with the orbits for which $t=0$, meaning that they
consist only of reflections and are supported on a single edge. Such
orbits exits only for even values of $n$ and their number is exactly
the number of edges, $E={V \choose 2}$; their repetition number cannot
be neglected, as it equals $r_{\p}=\nicefrac{n}{2}$.

Now, we refer to (\ref{eq:=000020form=000020factor=000020-=000020orbit=000020expansion})
for the form factor computation. We consider $n=2E\tau$, and next
we will eventually take the limit of increasing graphs. Namely, $E\rightarrow\infty$
(or equivalently $V\rightarrow\infty$ in our case), while fixing
the value of $\tau=\frac{2E}{n}$. When summing over pairs $\p,\q$
of periodic orbits in (\ref{eq:=000020form=000020factor=000020-=000020orbit=000020expansion})
we need to take into account only orbits with the same metric lengths,
i.e., $L_{\p}=L_{\q}$. Due to the complexity of this task, we do
not consider all such pairs, but only pairs in which either $\p=\q$
or that $\p$, $\q$ are the same up to a reversed direction (denoting
this by $\q=\hat{\p})$. This is the well-known diagonal approximation
and we will see that it successfully reproduces the leading term of
the GOE form factor. Another common approximation that we implement
here is to take $r_{\p}=r_{\q}=1$ for almost all orbits, apart from
the orbits for which we know their exact repetition number (these
are the orbits with $t=0$ and $r_{\p}=\nicefrac{n}{2}$, already
mentioned above). In the following computation we use the shorthand
notation $n:=2E\tau$, and return to the variable $\tau$ only at
its end. 

\begin{align}
K_{\textrm{diag}}(\tau)= & \frac{n^{2}}{2E}\sum_{\p\textrm{ has }n\textrm{ edges}}\frac{1}{r_{\p}^{2}}\left(A_{\p}A_{\p}^{*}+A_{\p}A_{\hat{\p}}^{*}\right)\label{eq:=000020Form-Factor-diagonal-complete-graphs}\\
= & \frac{n^{2}}{E}\sum_{\p\textrm{ has }n\textrm{ edges}}\frac{1}{r_{\p}^{2}}A_{\p}A_{\p}^{*}\nonumber \\
= & \frac{n^{2}}{E}\left[\frac{1}{2}\frac{1}{(n/2)^{2}}{V \choose 2}\left(1-\frac{2}{V-1}\right)^{2n}\right.\nonumber \\
 & \left.\quad\quad+\sum_{t=2}^{n}\frac{2}{n}{n \choose t}{V \choose 2}(V-2)^{t-2}\left(\frac{2}{V-1}\right)^{2t}\left(1-\frac{2}{V-1}\right)^{2(n-t)}\right]\nonumber \\
= & n^{2}\left[\frac{2}{n^{2}}\left(1-\frac{2}{V-1}\right)^{2n}+\sum_{t=2}^{n}\frac{2}{n}{n \choose t}(V-2)^{t-2}\left(\frac{2}{V-1}\right)^{2t}\left(1-\frac{2}{V-1}\right)^{2(n-t)}\right],\nonumber 
\end{align}
where in the above we used in the third line that for the first term
($t=0$) the orbits are such that $\p=\hat{\p}$ and hence there should
be introduced an additional factor of $\nicefrac{1}{2}$, and in the
last line we used $E={V \choose 2}$.

To continue the computation we note that a slight modification (starting
from $t=0$ rather than from $t=2$) of the sum above gives 
\begin{align*}
\sum_{t=0}^{n}\frac{2}{n}{n \choose t}(V-2)^{t-2} & \left(\frac{2}{V-1}\right)^{2t}\left(1-\frac{2}{V-1}\right)^{2(n-t)}\\
= & \frac{2}{n}\frac{1}{(V-2)^{2}}\sum_{t=0}^{n}{n \choose t}(V-2)^{t}\left(\frac{2}{V-1}\right)^{2t}\left(1-\frac{2}{V-1}\right)^{2(n-t)}\\
= & \frac{2}{n}\frac{1}{(V-2)^{2}}\left((V-2)\left(\frac{2}{V-1}\right)^{2}+\left(1-\frac{2}{V-1}\right)^{2}\right)^{n}\\
= & \frac{2}{n(V-2)^{2}}.
\end{align*}

We substitute this in (\ref{eq:=000020Form-Factor-diagonal-complete-graphs}),
while subtracting the additional terms ($t=0$ and $t=1$), and get 

\begin{align}
K_{\textrm{diag}}(\tau)= & n^{2}\left[\frac{2}{n^{2}}\left(1-\frac{2}{V-1}\right)^{2n}+\frac{2}{n(V-2)^{2}}\right.\label{eq:=000020Form-Factor-diagonal-complete-graphs-2}\\
 & \quad\quad\left.-\frac{2}{n}(V-2)^{-2}\left(1-\frac{2}{V-1}\right)^{2n}-2(V-2)^{-1}\left(\frac{2}{V-1}\right)^{2}\left(1-\frac{2}{V-1}\right)^{2(n-1)}\right]\nonumber \\
= & \frac{2n}{(V-2)^{2}}+2\left(1-\frac{2}{V-1}\right)^{2n}\left[1-\frac{n}{(V-2)^{2}}-\frac{4n^{2}}{(V-1)^{2}(V-2)}\left(1-\frac{2}{V-1}\right)^{-2}\right]\nonumber \\
 & \underset{E\rightarrow\infty}{\longrightarrow}2\tau,\nonumber 
\end{align}

where in the last line we used that $\tau=\frac{n}{2E}$ is fixed
and since $E={V \choose 2}$, we get that $\frac{n}{V^{2}}\rightarrow\tau$
as $V\rightarrow\infty$. To see this, one can observe that the second
term of the third line tends to $2\left(1-\frac{2}{V-1}\right)^{2\tau V^{2}}\left[1-\tau-4\tau^{2}V\right]$,
which goes to zero as $V\rightarrow\infty$. 

We now see that for small values of $\tau$, the calculation above
indeed reproduces the first term of the GOE form factor:
\begin{equation}
K_{\textrm{GOE}}(\tau)=\begin{cases}
2\tau-\tau\ln(1+2\tau) & \tau\leq1\\
2\tau-\tau\ln\left(\frac{2\tau+1}{2\tau-1}\right) & \tau\geq1.
\end{cases}\label{eq:=000020Form=000020Factor=000020-=000020Theoretic}
\end{equation}

We note that the computation above is valid also for the Neumann Kirchhoff
conditions. This is because the scattering coefficients of the Neumann-Kirchhoff
and of the asymptotic preferred orientation conditions are equal up
to sign (see (\ref{eq:=000020Even=000020degree=000020asymptotics})
and the text which follows it), but as we argued above, the sign is
canceled in the diagonal approximation.

In the computation above, before taking the limit, we may examine
the case $n=2$ (so that $\tau=\nicefrac{1}{E}$). Returning to (\ref{eq:=000020Form-Factor-diagonal-complete-graphs})
and taking only the first term there (since $t=0$ when $n=2$), we
get 
\begin{equation}
K_{\textrm{diag}}(\nicefrac{1}{E})=2\left(1-\frac{2}{V-1}\right)^{4}.\label{eq:=000020form=000020factor=000020for=000020n=00003D2}
\end{equation}
This explains the peak which may be observed in Figure~\ref{fig:=000020Form=000020factor=000020-=000020complete=000020graphs}
in the vicinity of zero and also explains why the peak for $K_{9}$
is higher than for $K_{7}$ (right and left parts of that figure).
The expression (\ref{eq:=000020form=000020factor=000020for=000020n=00003D2})
appeared already in \cite[eq. (31)]{SchSmi_pmb00}, where it is also
explained that such orbits, which are repetitions of 2-edge orbits.
are responsible for the odd-even staggering phenomena which ones obtains
numerically for small values of $\tau$.

\subsection{The peak at half the Heisenberg time}

\begin{figure}[h]
\begin{minipage}[t]{0.45\columnwidth}%
\includegraphics[width=1.12\textwidth]{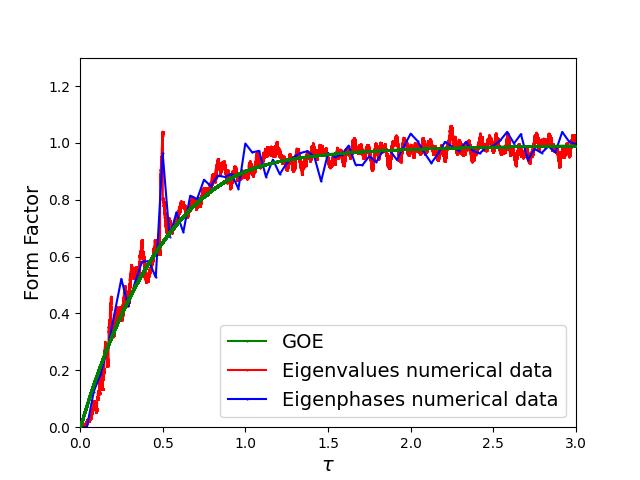}%
\end{minipage}$\quad\quad\quad$%
\begin{minipage}[t]{0.45\columnwidth}%
\includegraphics[width=1.12\textwidth]{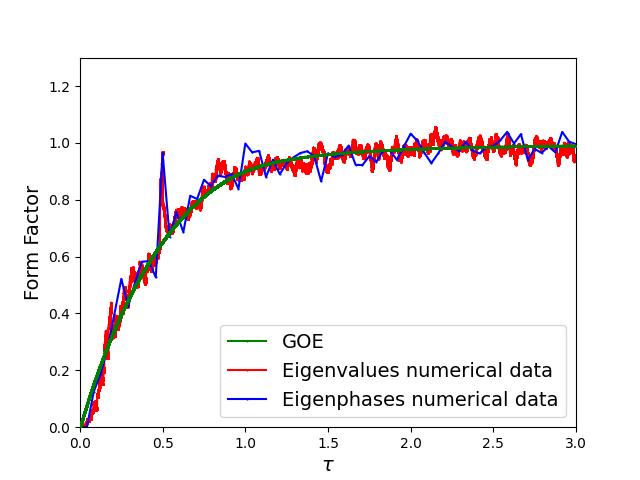}%
\end{minipage}

\caption{The form factor of the octahedron graph. Left: with preferred orientation
vertex conditions. Right: with Neumann-Kirchhoff vertex conditions.
The curves correspond to: numerics of eigenvalues (red), numerics
of the eigenphases (blue) and the theoretical GOE (green). \protect\label{fig:=000020Form=000020factor=000020of=000020Octahedron-2}}
\end{figure}

We return to the observation made in the introduction - when numerically
computing the form factor of the octahedron graph, we notice a clear
and substantial (though local) deviation from GOE (Figure~\ref{fig:=000020Form=000020factor=000020of=000020Octahedron-2},
Left). Specifically, there is a clear peak which appears at half the
Heisenberg time ($\tau=\nicefrac{1}{2}$). As far as we are aware
of, such phenomenon was never observed in form factors of chaotic
systems\footnote{To be accurate, a peak at $\tau=1/2$ may be spotted in \cite[fig. 2]{KotSch_phys01},
though there is no specific mention of this peak in the text.}. In what follows we explain the appearance of this peak and show
that there are additional setups (other graphs and other vertex conditions)
in which similar peaks appear. We start by considering a general graph
$\graph$ and upon need impose particular restrictions on the graph.
We refer to the periodic orbit expansion (\ref{eq:=000020form=000020factor=000020-=000020orbit=000020expansion})
and employ it to evaluate the value $K_{\SU}(\nicefrac{1}{2})$ and
to show that it is substantially higher than its GOE prediction, $K_{\textrm{GOE}}(\nicefrac{1}{2})=1-\nicefrac{1}{2}\ln(2)$.

Towards this we consider in the sum (\ref{eq:=000020form=000020factor=000020-=000020orbit=000020expansion})
only periodic orbits $\p,\q$ which consist of $\frac{2E}{2}=E$ edges.
Among such orbits we find the Eulerian cycles. An Eulerian cycle is
a closed path on the graph which visits every edge of the graph exactly
one time. For the sake of this definition we consider the undirected
graph, i.e., every undirected edge appears exactly once in an Eulerian
cycle. Therefore, the metric length of every Eulerian orbit equals
the total length of the graph, which means in particular that $L_{\p}=L_{\q}$
for every pair of two Eulerian orbits, $\p,\q$ and also that $r_{\p}=r_{\q}=1$.
Therefore, the contribution of Eulerian cycles to the sum in (\ref{eq:=000020form=000020factor=000020-=000020orbit=000020expansion})
is given by $\frac{1}{2}E\sum_{\p,\q\textrm{ Eulerian}}A_{\p}A_{\q}^{*}$.
To calculate the amplitude $A_{\p}$ for an Eulerian cycle $\p$,
first observe that (by definition) all scattering events in an Eulerian
cycle are transmissions. Therefore, by (\ref{eq:=000020Even=000020degree=000020asymptotics})
we get $S_{i,j}=-\frac{2}{d}(-1)^{i+j}$ for any scattering event
at any vertex. In addition, a vertex of even degree $d$ is being
visited exactly $\nicefrac{d}{2}$ times throughout an Eulerian cycle,
such that all edges connected to the vertex are eventually visited.
Hence, for such a vertex the product of all $\nicefrac{d}{2}$ scattering
amplitudes (which are all transmissions) gives
\[
\left(-\frac{2}{d}\right)^{d/2}\prod_{j=1}^{d}(-1)^{j}=\left(\frac{2}{d}\right)^{d/2},
\]
where we used that the degree $d$ is even and satisfies $d>2$, and
so both $\nicefrac{d}{2}$ and the sum $\sum_{j=1}^{d}j$ are even.
In order to simplify the following computations, we assume that the
graph $\graph$ is $d$-regular, i.e., all of its vertices are of
degree $d$. In such a case, multiplying all the scattering amplitudes
above for all the $V$ vertices we get that $A_{\p}=\left(\frac{2}{d}\right)^{Vd/2}=\left(\frac{2}{d}\right)^{E}$
for each Eulerian orbit $\p$. The relevant term in the form factor
is therefore 
\begin{equation}
\frac{1}{2}E\sum_{\p,\q\textrm{ Eulerian}}A_{\p}A_{\q}^{*}=\frac{1}{2}E\left(\frac{2}{d}\right)^{2E}\left(\Eucount(\graph)\right)^{2},\label{eq:=000020Eulerian=000020term=000020of=000020FF}
\end{equation}
where $\Eucount(\graph)$ is the number of the Eulerian cycles of
the graph $\graph$. To finish the computation we ought to count the
number of Eulerian cycles of a graph, $\Eucount(\graph)$. This is
a tedious task, which cannot be solved in a polynomial time for general
undirected graphs. We elaborate more about this problem and its possible
theoretical and practical resolutions in Subsection~\ref{subsec:=000020counting=000020Eulerian}
and in Appendix~\ref{sec:=000020algorithm=000020for=000020orbit=000020counting}.

If $\graph$ is taken to be the octahedron graph then $\Eucount(\textrm{Octahedron})=744$,
which we obtained by running two independent algorithms (as detailed
in Subsection~\ref{subsec:=000020counting=000020Eulerian} and in
appendix~\ref{sec:=000020algorithm=000020for=000020orbit=000020counting}).
To clarify, $\Eucount$ counts all the Eulerian orbits, where reversing
the direction of an orbit is considered as a new orbit, but a cyclic
shift of the orbit does not count as a new orbit. Substituting in
(\ref{eq:=000020Eulerian=000020term=000020of=000020FF}) gives
\begin{align*}
\frac{1}{2}E\sum_{\p,\q\textrm{ Eulerian}}A_{\p}A_{\q}^{*} & \approx0.198.
\end{align*}

The GOE form factor is $K_{\textrm{GOE}}(\nicefrac{1}{2})\approx0.6534$,
whereas the numerics gives $K_{\SU}(\nicefrac{1}{2})\approx0.9335\pm0.02$.
Hence, the contribution of the Eulerian cycles explains slightly more
than $70\%$ of this mismatch. This is satisfying since in general
the numerical values at $\tau=\frac{n}{2E}$ for even $n$ are higher
than the GOE prediction (and it is lower for odd values of $n$),
which is due to the staggering phenomenon mentioned at the end of
the previous subsection. 

We note that the octahedron has many periodic orbits which consist
of $12$ edges; there are more than $1.4\cdot10^{6}$ of those. We
get this number by computing $\frac{1}{12}\tr{C^{12}}=1398784$, where
$C$ is the $6\times6$ adjacency (connectivity) matrix of the octahedron
graph. But, we take in consideration that computing via the trace
of $C^{12}$ underestimates the orbits $\p$ with repetition number
$r_{\p}>1$ (for those orbits we should divide the trace by $12/r_{p}$
and not by the global factor $12$ as we do). Overall, we see that
the number of Eulerian cycles is less than $0.05\%$ from all the
periodic orbits of size $12$ (of the octahedron). Their substantial
effect on the form factor is not because of their number, but it is
thanks to the constructive interference between \emph{every} two such
periodic orbits. Namely, for the Eulerian cycles all the terms in
the sum $\sum_{\p,\q\textrm{ Eulerian}}A_{\p}A_{\q}^{*}$ are of positive
sign and hence their contribution is constructive. This explains the
dominance of the Eulerian cycles in this case and the dominance of
the observed peek at $\tau=1/2$.

We complement this computation by mentioning that changing the vertex
conditions from preferred orientation to Neumann-Kirchhoff does not
affect the result and the existence of the peak (see Figure~\ref{fig:=000020Form=000020factor=000020of=000020Octahedron-2},
Right). Indeed, the transmission amplitude of a scattering event for
Neumann-Kirchhoff conditions is $2/d$, which equals\footnote{Up to sign which is always canceled in our case, as was already mentioned.}
the transmission coefficient in the preferred orientation conditions.
Therefore the contribution of the Eulerian cycles computed in (\ref{eq:=000020Eulerian=000020term=000020of=000020FF})
is exactly the same if preferred orientation conditions are changed
into Neumann-Kirchhoff.

\begin{figure}[h!]
\begin{minipage}[t]{0.45\columnwidth}%
\includegraphics[scale=0.5]{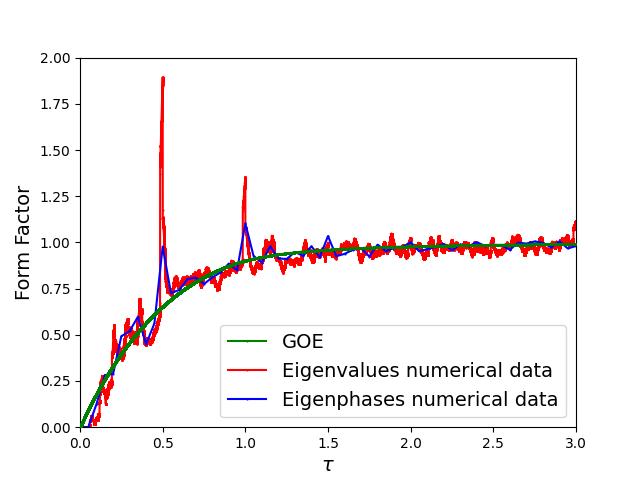}%
\end{minipage}\hfill{}%
\begin{minipage}[t]{0.45\columnwidth}%
\includegraphics[scale=0.5]{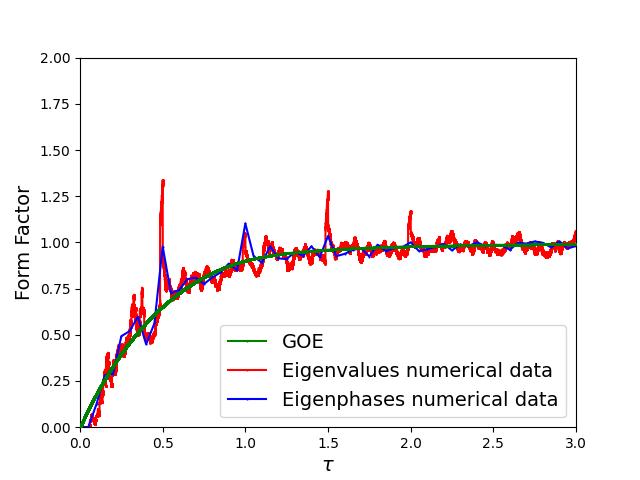}%
\end{minipage}\caption{The form factor of the $K_{5}$ graph. Numerics of eigenvalues (red),
numerics of the eigenphases (blue) and the theoretical GOE (green).
Left: with preferred orientation vertex conditions. Right: with Neumann-Kirchhoff
vertex conditions.\protect\label{fig:=000020Form=000020factors=000020with=000020peaks=000020-=000020more}}
\end{figure}

Experimentally, we have observed the appearance of such a peak at
$\tau=\nicefrac{1}{2}$ also for the complete graph $K_{5}$ (see
Figure~\ref{fig:=000020Form=000020factors=000020with=000020peaks=000020-=000020more}).
Repeating the computation above for $K_{5}$ gives 
\begin{align}
\frac{1}{2}E\sum_{\p,\q\textrm{ Eulerian}}A_{\p}A_{\q}^{*}= & \frac{1}{2}E\left(\frac{2}{d}\right)^{2E}\left(\Eucount(K_{5})\right)^{2}\approx0.332,\label{eq:=000020K_5=000020Eulerian=000020contribution}
\end{align}
where we used that $K_{5}$ is a $d$-regular graph with $d=4$, $E=10$,
and $\Eucount(K_{5})=264$. The computation in (\ref{eq:=000020K_5=000020Eulerian=000020contribution})
agrees with the difference between the numerical (eigenphases) and
the theoretical value, which is about $0.324$.

Furthermore, we observe in the form factor for $K_{5}$ similar peaks
at multiples of $\tau=\nicefrac{1}{2}$. A first attempt to explain
those might be via orbits which are concatenation of a few Eulerian
cycles (not necessarily a repetition of the same Eulerian cycle, but
rather combining a few). But, such computations do not yield satisfactory
values; the problem of providing a complete explanation for these
particular peaks is still open.

We end this section by pointing out that the peak at $\tau=\nicefrac{1}{2}$
(and sometimes at its multiple) does not appear for every graph. For
example, we do not see such peaks for graphs which are not Eulerian
(i.e., they have no Eulerian cycles). A graph is Eulerian if and only
if all of its vertices are of even degree; so complete graphs of the
form $K_{2n}$ are not Eulerian and indeed do not show these peaks.
Nevertheless, even in the form factor of Eulerian graphs we do not
necessarily see the aforementioned peaks. To give an example, this
is the case with $K_{7}$ and $K_{9}$, as may be seen in Figure~\ref{fig:=000020Form=000020factor=000020-=000020complete=000020graphs}.
If we repeat the computation in (\ref{eq:=000020K_5=000020Eulerian=000020contribution})
for $K_{7}$, substituting $E=21$, $d=6$ and $\Eucount(K_{7})=129976320$,
we get that Eulerian cycles contribute $0.00162$ to the form factor
at $\tau=\nicefrac{1}{2}$. This is indeed negligible comparing to
the GOE value, $K_{\textrm{GOE}}(\nicefrac{1}{2})\approx0.6534$,
so no peak is seen in this case.

Similarly, repeating the computation in (\ref{eq:=000020K_5=000020Eulerian=000020contribution})
for $K_{9}$, substituting $E=36$, $d=8$ and $\Eucount(K_{9})=911520057021235200$,
we get that Eulerian cycles contribute $6.7\cdot10^{-7}$ to the form
factor at $\tau=\nicefrac{1}{2}$, which is again negligible comparing
to the GOE value.

We end by returning to the asymptotic computation of complete graphs
from the previous subsection and checking whether the peak at $\tau=\nicefrac{1}{2}$
appears there. The asymptotic formula for the number of Eulerian cycles
in complete graphs $K_{n}$ (where $n$ is odd and tends to infinity)
is known to be (see \cite{McKayRob_cpc98}), 
\[
\Eucount(K_{n})=2^{(n+1)/2}\pi^{1/2}\rme^{-n^{12}/2+11/12}n^{(n-2)(n+1)/2}(1+O(n^{-1/2+\epsilon})).
\]

Using this asymptotics in (\ref{eq:=000020K_5=000020Eulerian=000020contribution})
with $E={n \choose 2}$ and $d=n-1$ gives that the contribution to
the form factor is of order
\[
2^{(2n^{2}-n-3)/2}\pi^{1/2}\rme^{-n^{12}/2+11/12}n^{-(n-3)(n+2)/2}(1+O(n^{-1/2+\epsilon})),
\]
which converges to $0$, as $n$ tends to infinity. 

\subsection{The problem of counting Eulerian cycles\protect\label{subsec:=000020counting=000020Eulerian}}

The general problem of counting the exact number of Eulerian cycles
of an \emph{undirected} graph belongs to the class of $\#p$-complete
problems \cite{BriWin_2005}. This means in particular that there
is no known polynomial time algorithm which solves it. Furthermore,
if such polynomial solution is found it would imply that all $\#p$
problems\footnote{These are analogous to $NP$ problems, but instead of just finding
whether a solution to the problem exists, one is required to counting
the number of all the possible solutions.} could also be solved in polynomial time, which would lead to significant
implications in computational complexity theory. We do not aim here
to cover the background on this computational problem and merely mention
a very recent preprint which offers a new approach and reviews the
existing literature \cite{Luo_arxiv25}. The problem of counting periodic
orbits and even estimating their number or providing asymptotics plays
an important role in quantum chaos, as it allows a better control
when using trace formulae for the study of spectral statistics. The
orbit counting problem is relevant for graphs, as well as for the
case of symbolic dynamics. In both cases one wishes to cluster the
periodic orbits according to the relevant problem. For example, in
the case of metric graphs one would like to cluster all orbits According
to their support. Such a cluster would contain all orbits, which share
the number of times they transverse each of the graph edges. Using
this approach and offering numerical as well as analytical solutions
to it was done in \cite{KotSch_phys01,EchHarHud_grcomb25,Tan_jpa00,Ber_incol06,SchSmi_pmb00,SchSmi_ejc01,GavSmi_jpa07}
for specific families of metric graphs and in \cite{GutOsi_nonlin13,GutOsi_jsp13}
for general systems, by considering symbolic dynamics.

We mention here two different algorithms which we have used in order
to exactly count the number of Eulerian cycles in the graphs considered
in this paper.

First, we note that in the case of a directed graph there is a well-known
polynomial time algorithm to count the number of Eulerian cycles.
This algorithm is based on the BEST theorem (after de Bruijn, van
Aardenne-Ehrenfest, Smith, and Tutte) which provides an explicit formula
for the number of Eulerian cycles in terms of the number of spanning
trees of a graph \cite{EhrBru_sis1951,SmiTut_amm1941}. One may apply
such an algorithm straightforwardly to our problem (i.e., counting
for non-directed graphs): enumerate over all possible assignments
of directions of edges and for each perform the BEST algorithm. This
is clearly exponential in the number of graph edges due to the enumeration,
but may still be performed for small enough graphs. Indeed, we used
this for the computation of the $\Eucount$ values given in the previous
subsection.

The other algorithm we have used is new, to the best of our knowledge,
and we describe it in detail in Appendix~\ref{sec:=000020algorithm=000020for=000020orbit=000020counting}.
Both algorithms are of similar order of complexity, but the algorithm
provided in the appendix may be used to count other families of periodic
orbits beyond Eulerian cycles and also to evaluate their contribution
to the spectral statistics.

\section*{Acknowledgments}

We greatly benefited from interesting discussions with Gregory Berkolaiko
and Uzy Smilansky on random matrix theory in the context of quantum
graphs. We had stimulating discussions with Sven Gnutzmann, Boris
Gutkin and Holger Schanz about periodic orbit expansions. We thank
all of them for enriching our knowledge on these topics. Special thanks
are due to Jon Harrison for his thoughtful feedback on this manuscript.

R.B. was supported by the Israel Science Foundation (ISF Grant No.
844/19) and by the Binational Foundation Grant (grant no. 2016281).
P.E. was partially supported by the European Union’s Horizon 2020
research and innovation programme under the Marie Skłodowska-Curie
grant agreement No 873071.

\appendix

\section{An algorithm for counting periodic orbits\protect\label{sec:=000020algorithm=000020for=000020orbit=000020counting}}

We present here a procedure for counting families periodic orbits.
We applied it numerically or the counting of Eulerian cycles which
play important role in the current work. Nevertheless, this scheme
may be applied for counting various types of orbits.

\subsection{The adjacency matrix for orbit counting}

Let $\graph$ be a directed graph with the vertex set $\vtcs$ and
the directed edge set $\edges$. We consider here edge sets with the
property $e\in\edges\Leftrightarrow\hat{e}\in\edges$, where $\hat{e}$
denotes the reverse direction of an edge $e$. As a starting point
for counting periodic orbits we use the graph adjacency matrix $\adj(\graph)$,
which is a $V\times V$ matrix such that 
\begin{equation}
\left[\adj(\graph)\right]_{i,j}=\begin{cases}
1 & (j,i)\textrm{ is a directed edge}\\
0 & \textrm{otherwise}
\end{cases},\label{eq:=000020adjacency=000020of=000020entire=000020graph}
\end{equation}
where we assume that two vertices may be connected by at most a single
edge. Because every directed edge also appears with its reverse direction,
we get that the matrix $\adj$ is symmetric, $\adj_{i,j}=\adj_{j,i}$.
Observe that the number of periodic orbits of length $n$ of the graph
is characterized by $\tr{\adj^{n}}$. Here, one should be careful
because some orbits are counted with multiplicity. Specifically, the
periodic orbits are counted in $\tr{\adj^{n}}$ with some factor which
is due to taking cyclic permutations. Namely, a certain orbit $\left(v_{1},v_{2},\ldots,v_{n},v_{1}\right)$
with all vertices different $v_{i}\neq v_{j}$ will be counted $n$
times in $\tr{\adj^{n}}$, since all its cyclic permutations $\left(v_{k},v_{k+1},\ldots,v_{n},v_{1},v_{2},\ldots,v_{k}\right)$,
for $1\leq k\leq n$, will be counted. There are particular orbits
which are repetitions of shorter orbits, such as $\p=(v_{1},v_{2},\ldots,v_{m},v_{1},v_{2},\ldots,v_{m},\ldots,$\allowbreak$v_{1},v_{2},\ldots,v_{m})$.
Here we assume as before that $v_{i}\neq v_{j}$ for $1\leq i,j\leq m$
and that repeated vertices are indicated explicitly. In such a case
we get that $m\mid n$ (i.e., $m$ divides $n$) and further denote
$r_{\p}:=\frac{n}{m}$, which is called the repetition number of the
orbit $\p$. Using this notation, we see that each periodic orbit
of length $n$ is counted exactly $n/r_{\p}$ times in $\tr{\adj^{n}}$.
We summarize this by writing
\begin{equation}
\tr{\adj^{n}}=\sum_{\p\in\PO n}\frac{n}{r_{\p}}=n\sum_{\p\in\PO n}\frac{1}{r_{p}},\label{eq:=000020trace=000020of=000020A^n=000020via=000020periodic=000020orbits}
\end{equation}
where $\PO n$ is the set of periodic orbits of length $n,$ and $r_{\p}$
is the repetition number of the orbit $\p$.

\subsection{counting orbits of subgraphs}

Next, we wish to count the orbits of subgraphs of $\graph$. To do
so, we repeat exactly the same arguments of the previous subsection,
but for a modified adjacency matrix. Explicitly, let $\widetilde{\graph}$
be a subgraph of $\graph$ having the same vertex set $\vtcs$ as
$\graph$, but only a subset of the directed edges $\widetilde{\edges}\subset\edges$
(for example, in this subset it might happen that a certain edge appears
but its reverse does not). We construct the adjacency matrix $\adj(\widetilde{\graph})$
of the subgraph as a $V\times V$ matrix such that 
\[
\left[\adj(\widetilde{\graph})\right]_{i,j}=\begin{cases}
1 & (j,i)\in\widetilde{\edges}\\
0 & \textrm{otherwise}
\end{cases}.
\]
 We still have that the analogue of (\ref{eq:=000020trace=000020of=000020A^n=000020via=000020periodic=000020orbits})
holds, but obviously in the sum on the right hand side only periodic
orbits which are supported on the edges set $\widetilde{\edges}$
are taken into account.

\subsection{Introducing vectors of counts}

We introduce the following notation for the power set of the edge
set, $\bonds:=\ZZ^{\left|\edges\right|}$. With this notation there
is a bijection between subsets $\widetilde{\edges}\subset\edges$
and $\bvec\in\bonds$, and we may think of the latter as sequences
of bits of length $\left|\edges\right|$. From now on, we fix $n$
(the length of the periodic orbits we are counting) and denote
\begin{equation}
\widetilde{N}_{\bvec}:=\tr{\adj(\widetilde{\graph})^{n}},\label{eq:=000020vec_tilde_as_trace}
\end{equation}
where $\bvec\in\bonds$ corresponds to the edge subset $\widetilde{\edges}$
of the subgraph $\widetilde{\Gamma}$. We note that $\widetilde{N}_{\bvec}$
counts the periodic orbits whose support is $\bvec$ or a subset of
\textbf{$\bvec$} (where $\bvec$ is understood from now on as an
edge subset).

We denote by $\supp(\p)$ the subset of edges which the periodic orbit
$\p$ contains, and allow ourselves to write $\supp(\p)\in\bonds$,
thanks to the bijection between $\bonds$ and edge subsets. It might
be that a certain edge appears more than once in $\p$, but such information
is not reflected in the notation $\supp(\p)$.

With that notation we write (following (\ref{eq:=000020trace=000020of=000020A^n=000020via=000020periodic=000020orbits})),
\begin{equation}
\widetilde{N}_{\bvec}=n\sum_{\supp(p)\subseteq\bvec}\frac{1}{r_{\p}},\label{eq:=000020vec_tilde}
\end{equation}
and 
\begin{equation}
N_{\bvec}:=n\sum_{\supp(p)=\bvec}\frac{1}{r_{\p}},\label{eq:=000020vec_exact}
\end{equation}
where we sum only over $\p\in\PO n\textrm{ },$and noting that the
difference between $\widetilde{N}_{\bvec}$ and $N_{\bvec}$ is whether
the considered periodic orbits are supported exactly on $\bvec$ or
on some subset of it.

We consider the vectors $\left(\widetilde{N}_{\bvec}\right)_{\bvec\in\bonds}$
and $\left(N_{\bvec}\right)_{\bvec\in\bonds}$ as vectors of lengths
$2^{\left|\edges\right|}$. Namely, $\widetilde{N},N\in\N_{0}^{2^{\left|\edges\right|}}$,
where $\N_{0}:=\N\cup\{0\}$.

The vector $\left(N_{\bvec}\right)_{\bvec\in\bonds}$ may be used
to count various families of periodic orbits. As an example, we demonstrate
how it can be used to count Eulerian cycles. We start by observing
that there are necessary conditions which a certain $\bvec\in\bonds$
must fulfill in order to be the support of an Eulerian cycle. Namely,
$\bvec$ should correspond to a subset $\widetilde{\edges}\subset\edges$
which contains exactly $\left|\edges\right|/2$ edges and satisfies
$e\in\widetilde{\edges}\Leftrightarrow\hat{e}\notin\widetilde{\edges}$
(but this is not a sufficient condition). We denote the set of these
admissible $\bvec$ values by $\Eubonds$. Using this, the number
of Eulerian cycles of the graph $\graph$ is 
\begin{equation}
\Eucount(\graph)=\frac{2}{\left|\edges\right|}\sum_{\bvec\in\Eubonds}N_{\bvec},\label{eq:=000020Euler=000020count=000020in=000020terms=000020of=000020vec}
\end{equation}
where we used (\ref{eq:=000020vec_exact}) and that $n=\left|\edges\right|/2$
and that $r_{\p}=1$ for an Eulerian cycle $\p$.

Next, we describe how to express the vector $\left(N_{\bvec}\right)_{\bvec\in\bonds}$
in terms of the vector $\left(\widetilde{N}_{\bvec}\right)_{\bvec\in\bonds}$
. The latter vector was already explicitly expressed in (\ref{eq:=000020vec_tilde_as_trace}).

\subsection{Transform from $\widetilde{N}_{\protect\bvec}$ to $N_{\protect\bvec}$}

We keep in mind that $1\leq n\leq\left|\edges\right|$ is fixed throughout
the section (for Eulerian cycles one takes $n=\left|\edges\right|/2$
, but we continue describing the scheme for an arbitrary value of
$n$). By definition and from (\ref{eq:=000020vec_tilde}) and (\ref{eq:=000020vec_exact}),
one observes that the values $\left(N_{\bvec}\right)_{\bvec\in\bonds}$
may be expressed in terms of $\left(\widetilde{N}_{\bvec}\right)_{\bvec\in\bonds}$
by using the inclusion--exclusion principle. To write this, we introduce
a few additional notations. Denote by $\pop(\bvec)$ the number of
$1$'s which $\bvec$ contains. Furthermore, if $\widetilde{\bvec},\widehat{\bvec}\in\bonds$
correspond to some edge subsets $\widetilde{\edges},\widehat{\edges}\subset\edges$
then we denote $\widetilde{\bvec}\subset\widehat{\bvec}$ iff $\widetilde{\edges}\subset\widehat{\edges}$.
We further emphasize this by writing that $\widetilde{\bvec}\subset\widehat{\bvec}$
iff $\widetilde{\bvec}\&\widehat{\bvec}=\widetilde{\bvec}$, where
notation $\&$ is the pair-wise 'AND' operation between bits. Given
$\bvec\in\bonds$ we denote for all $k\in\N$, 
\[
\bbonds k:=\set{\btvec\in\bonds}{\btvec\subset\bvec\textrm{ and }\pop(\btvec)=\pop(\bvec)-k}.
\]
Using these notation we employ the inclusion--exclusion principle
and get 
\begin{equation}
N_{\bvec}=\widetilde{N}_{\bvec}-\sum_{\btvec\in\bbonds 1}\widetilde{N}_{\btvec}+\sum_{\btvec\in\bbonds 2}\widetilde{N}_{\btvec}-\ldots+(-1)^{n-1}\sum_{\btvec\in\bbonds{n-1}}\widetilde{N}_{\btvec}+(-1)^{n}\sum_{\btvec\in\bbonds n}\widetilde{N}_{\btvec},\label{eq:=000020vec_exact=000020as=000020linear=000020combination=000020of=000020vec_tilde}
\end{equation}
where we use that $\pop(\bvec)\leq n$. Note also that $\bbonds k=\emptyset$
for $k>\pop(\bvec)$ and $\bbonds{\pop(\bvec)}=\left\{ \left(0,\ldots,0\right)\right\} $
and $\widetilde{N}_{\left(0,\ldots,0\right)}=0$; which means in particular
that the last summand in (\ref{eq:=000020vec_exact=000020as=000020linear=000020combination=000020of=000020vec_tilde})
always vanishes $(-1)^{n}\sum_{\btvec\in\bbonds n}\widetilde{N}_{\btvec}=0$.

Let us consider $\left(N_{\bvec}\right)_{\bvec\in\bonds}$ and $\left(\widetilde{N}_{\btvec}\right)_{\btvec\in\bonds}$
as vectors in $\R^{\bonds}$, recalling that $\left|\bonds\right|=2^{\left|\edges\right|}$.
Equation (\ref{eq:=000020vec_exact=000020as=000020linear=000020combination=000020of=000020vec_tilde})
presents the vector $\left(N_{\bvec}\right)_{\bvec\in\bonds}$ as
a linear transform of the vector $\left(\widetilde{N}_{\btvec}\right)_{\btvec\in\bonds}$.
Denoting the matrix representing this linear transformation by $\mat{\left|\edges\right|}$,
we see that its entries are $0$, $\pm1$ and furthermore it has a
recursive representation as 
\begin{align}
\mat 1 & :=\left(\begin{array}{cc}
1 & 0\\
-1 & 1
\end{array}\right),\nonumber \\
\mat n & =\mat 1\otimes\mat{n-1},\quad\forall n>1,\label{eq:=000020Kronecker=000020Product}
\end{align}
where $\otimes$ denotes the Kronecker product of matrices. This transform
is also known as the arithmetic transform, which is closely related
to the well-known Walsh-Hadamard transform \cite{Arndt_2011}.

This special form of the matrices allows an efficient algorithm which
multiplies the matrix $\mat{\left|\edges\right|}$ with the vector
$\left(\widetilde{N}_{\btvec}\right)_{\btvec\in\bonds}$ to obtain
the vector $\left(N_{\bvec}\right)_{\bvec\in\bonds}$. This algorithm
uses a method similar to the fast Walsh-Hadamard. Its complexity is
$\left|\bonds\right|\log\left(\left|\bonds\right|\right)=\left|\edges\right|2^{\left|\edges\right|}$
instead of the usual complexity of multiplying a matrix by a vector
(which is $\left|\bonds\right|^{2}=2^{2\left|\edges\right|}$ in our
case).

\subsection{Algorithmic summary and complexity}

We summarize the steps of the algorithm described above for computing
$\Eucount(\graph)$ and its overall complexity:
\begin{enumerate}
\item Preparing the vector $\left(\widetilde{N}_{\bvec}\right)_{\bvec\in\bonds}$
, as in (\ref{eq:=000020vec_tilde_as_trace}). We need to perform
$\left|\bonds\right|=2^{\left|\edges\right|}$ times the computation
$\tr{\adj(\widetilde{\graph})^{\left|\edges\right|/2}}$ (with a different
sub-graph $\widetilde{\graph}$ every time). Overall the complexity
is $O\left(2^{\left|\edges\right|}\left|\vtcs\right|^{k}\log\left|\edges\right|\right)$,
with $2.37\lesssim k\leq3$ ($k$ depends on the multiplication algorithm
we choose, see e.g. \cite{AmbFilLeGa_proc15}).
\item Transforming $\left(\widetilde{N}_{\bvec}\right)_{\bvec\in\bonds}$
into $\left(N_{\bvec}\right)_{\bvec\in\bonds}$ using the arithmetic
transform given by (\ref{eq:=000020Kronecker=000020Product}). The
complexity is $O\left(\left|\bonds\right|\log\left(\left|\bonds\right|\right)\right)=O\left(\left|\edges\right|2^{\left|\edges\right|}\right)$.
\item Summing the relevant entries of $\left(N_{\bvec}\right)_{\bvec\in\bonds}$,
see (\ref{eq:=000020Euler=000020count=000020in=000020terms=000020of=000020vec}).
Overall it means to sum $\left|\Eubonds\right|=2^{\nicefrac{\left|\edges\right|}{2}}$
entries. 
\end{enumerate}
The total complexity is hence $O\left(2^{\left|\edges\right|}\left(\left|\edges\right|+\left|\vtcs\right|^{k}\log\left|\edges\right|\right)\right)$
with $2.37\lesssim k\leq3$ (depending on the multiplication algorithm).

\subsection{Concluding remarks and variations on the algorithm}

A first variation on the scheme described above would be to replace
the adjacency matrix $A$ in (\ref{eq:=000020adjacency=000020of=000020entire=000020graph})
by the overall $\left|\edges\right|\times\left|\edges\right|$ scattering
matrix, $\mathbf{S}$ or some modifications of it. For example, we
may use an $\left|\edges\right|\times\left|\edges\right|$ matrix
$\mathbf{M}$, defined by $\mathbf{M}_{i,j}=\text{\ensuremath{\mathbf{S}_{i,j}^{2}} }$
to evaluate the diagonal approximation as in (\ref{eq:=000020Form-Factor-diagonal-complete-graphs}).
Specifically, one may take the analogue of (\ref{eq:=000020vec_tilde_as_trace})
to define $\widetilde{N}_{\bvec}:=\tr{\mathbf{M}(\widetilde{\graph})^{n}}$.
Transforming $\widetilde{N}_{\bvec}$ (by the same arithmetic transform)
gives the vector $N_{\bvec}$, the sum of whose entries is the corresponding
value of the diagonal approximation.

Another variation of the algorithm is relevant for its implementation.
Rather then fixing a single value of $n$ and raising all matrices
to this power $n$, one may take several such $n$ values and correspondingly
prepare several vectors $\left(\widetilde{N}_{\bvec}\right)_{\bvec\in\bonds}$
(one for each power $n$). Then the arithmetic transform to turn them
into $\left(N_{\bvec}\right)_{\bvec\in\bonds}$ may be done in parallel,
thus making the computation more efficient.

\bibliographystyle{alpha}
\bibliography{GlobalBib_250731}

\newcommand{\etalchar}[1]{$^{#1}$}
\begin{thebibliography}{HBPan{\etalchar{+}}04}

\bibitem[AFLG15]{AmbFilLeGa_proc15}
A.~Ambainis, Y.~Filmus, and F.~Le~Gall.
\newblock Fast matrix multiplication: Limitations of the laser method.
\newblock In {\em Proceedings of the 47th Annual ACM Symposium on Theory of
  Computing (STOC)}, pages 585--593, 2015.

\bibitem[AG15]{AkiGut_jpa15}
M.~Akila and B.~Gutkin.
\newblock Spectral statistics of nearly unidirectional quantum graphs.
\newblock {\em J. Phys. A}, 48(34):345101, 21, 2015.

\bibitem[AG19]{Akiut_jpa19}
M.~Akila and B.~Gutkin.
\newblock G{SE} spectra in uni-directional quantum systems.
\newblock {\em J. Phys. A}, 52(23):235201, 8, 2019.

\bibitem[Arn11]{Arndt_2011}
J.~Arndt.
\newblock {\em The Walsh transform and its relatives}, pages 459--496.
\newblock Springer Berlin Heidelberg, Berlin, Heidelberg, 2011.

\bibitem[BE09]{BolEnd_ahp09}
J.~Bolte and S.~Endres.
\newblock The trace formula for quantum graphs with general self adjoint
  boundary conditions.
\newblock {\em Ann. Henri Poincar\'e}, 10(1):189--223, 2009.

\bibitem[Ber00]{Berkolaiko_thesis00}
G.~Berkolaiko.
\newblock {\em Quantum star graphs and related systems}.
\newblock Phd thesis, University of Bristol, 2000.

\bibitem[Ber04]{Ber_wrm04}
G.~Berkolaiko.
\newblock Form factor for large quantum graphs: evaluating orbits with time
  reversal.
\newblock {\em Waves Random Media}, 14(1):S7--S27, 2004.
\newblock Special section on quantum graphs.

\bibitem[Ber06]{Ber_incol06}
G.~Berkolaiko.
\newblock Form factor expansion for large graphs: a diagrammatic approach.
\newblock In {\em Quantum graphs and their applications}, volume 415 of {\em
  Contemp. Math.}, pages 35--49. Amer. Math. Soc., Providence, RI, 2006.

\bibitem[Ber17]{Berkolaiko_qg-intro17}
G.~Berkolaiko.
\newblock An elementary introduction to quantum graphs.
\newblock In {\em Geometric and computational spectral theory}, volume 700 of
  {\em Contemp. Math.}, pages 41--72. Amer. Math. Soc., Providence, RI, 2017.

\bibitem[BG00]{BarGas_jsp00}
F.~Barra and P.~Gaspard.
\newblock On the level spacing distribution in quantum graphs.
\newblock {\em J. Statist. Phys.}, 101(1--2):283--319, 2000.

\bibitem[BG18]{BanGnu_qg-exerices18}
R.~Band and S.~Gnutzmann.
\newblock Quantum graphs via exercises.
\newblock In {\em Spectral theory and applications}, volume 720 of {\em
  Contemp. Math.}, pages 187--203. Amer. Math. Soc., Providence, RI, 2018.

\bibitem[BGS84]{BohGiaSch_prl84}
O.~Bohigas, M.~J. Giannoni, and C.~Schmit.
\newblock Characterization of chaotic quantum spectra and universality of level
  fluctuation laws.
\newblock {\em Phys. Rev. Lett.}, 52(1):1--4, 1984.

\bibitem[BH03a]{BolHar_jpa03a}
J~Bolte and J~Harrison.
\newblock Spectral statistics for the {D}irac operator on graphs.
\newblock {\em J. Phys. A}, 36(11):2747--2769, 2003.

\bibitem[BH03b]{BolHar_jpa03}
J.~Bolte and J.~Harrison.
\newblock Ths spin contribution to the form factor of quantum graphs.
\newblock {\em J. Phys. A}, 36(27):L433--L440, 2003.

\bibitem[BH06]{BolHar_incol06}
J.~Bolte and J.~Harrison.
\newblock The spectral form factor for quantum graphs with spin-orbit coupling.
\newblock In {\em Quantum graphs and their applications}, volume 415 of {\em
  Contemp. Math.}, pages 51--64. Amer. Math. Soc., Providence, RI, 2006.

\bibitem[BHJ12]{BanHarJoy_jpa12}
R.~Band, J.~M. Harrison, and C.~H. Joyner.
\newblock Finite pseudo orbit expansions for spectral quantities of quantum
  graphs.
\newblock {\em J. Phys. A}, 45(32):325204, 19, 2012.

\bibitem[BHS19]{BanHarSep_jmaa19}
R.~Band, J.~M. Harrison, and M.~Sepanski.
\newblock Lyndon word decompositions and pseudo orbits on {$q$}-nary graphs.
\newblock {\em J. Math. Anal. Appl.}, 470(1):135--144, 2019.

\bibitem[BK99]{BerKea_jpa99}
G.~Berkolaiko and J.~P. Keating.
\newblock Two-point spectral correlations for star graphs.
\newblock {\em J. Phys. A}, 32(45):7827--7841, 1999.

\bibitem[BK13]{BerKuc_graphs}
G.~Berkolaiko and P.~Kuchment.
\newblock {\em Introduction to Quantum Graphs}, volume 186 of {\em Math. Surv.
  and Mon.}
\newblock AMS, 2013.

\bibitem[BSW02]{BerSchWhi_prl02}
G.~Berkolaiko, H.~Schanz, and R.~S. Whitney.
\newblock Leading off-diagonal correction to the form factor of large graphs.
\newblock {\em Phys. Rev. Lett.}, 88:104101, 2002.

\bibitem[BSW03]{BerSchWhi_jpa03}
G.~Berkolaiko, H.~Schanz, and R.~S. Whitney.
\newblock Form factor for a family of quantum graphs: an expansion to third
  order.
\newblock {\em J. Phys. A}, 36(31):8373--8392, 2003.

\bibitem[BW05]{BriWin_2005}
G.~R. Brightwell and P.~Winkler.
\newblock Counting eulerian circuits is \#p-complete.
\newblock In {\em ALENEX/ANALCO}, 2005.

\bibitem[BW10]{BerWin_tams10}
G.~Berkolaiko and B.~Winn.
\newblock Relationship between scattering matrix and spectrum of quantum
  graphs.
\newblock {\em Trans. Amer. Math. Soc.}, 362(12):6261--6277, 2010.

\bibitem[CGK{\etalchar{+}}]{CheGluKohGuhDie_arXiv25}
J.~Che, N.~Gluth, S.~Köhnes, T.~Guhr, and B.~Dietz.
\newblock Experimental study of the distributions of off-diagonal
  scattering-matrix elements of quantum graphs with symplectic symmetry.
\newblock {\em arXiv:2505.09573}.

\bibitem[Dav79]{Davis1979}
P.~J. Davis.
\newblock {\em Circulant matrices}.
\newblock A Wiley-Interscience Publication. John Wiley \& Sons, New
  York-Chichester-Brisbane, 1979.
\newblock Pure and Applied Mathematics.

\bibitem[DKM{\etalchar{+}}24]{DieKlaMasMisRicSkiWun_pre24}
B.~Dietz, T.~Klaus, M.~Masi, M.~Miski-Oglu, A.~Richter, T.~Skipa, and
  M.~Wunderle.
\newblock Closed and open superconducting microwave waveguide networks as a
  model for quantum graphs.
\newblock {\em Phys. Rev. E}, 109:034201, Mar 2024.

\bibitem[DYB{\etalchar{+}}17]{DieYunBiaBauLawSir_pre17}
B.~Dietz, V.~Yunko, M.~Bia\l{}ous, S.~Bauch, M.~\L{}awniczak, and L.~Sirko.
\newblock Nonuniversality in the spectral properties of time-reversal-invariant
  microwave networks and quantum graphs.
\newblock {\em Phys. Rev. E}, 95:052202, May 2017.

\bibitem[EHH25]{EchHarHud_grcomb25}
I.~Echols, J.~Harrison, and T.~Hudgins.
\newblock Periodic orbits on 2-regular circulant digraphs.
\newblock {\em Graphs Combin.}, 41(3):Paper No. 67, 17, 2025.

\bibitem[ET18]{Exner2018a}
P.~Exner and M.~Tater.
\newblock Quantum graphs with vertices of a preferred orientation.
\newblock {\em Physics Letters A}, 382(5):283--287, 2018.

\bibitem[ET21]{Exner2021}
P.~Exner and M.~Tater.
\newblock Quantum graphs: self-adjoint, and yet exhibiting a nontrivial
  {$\mathcal{PT}$}-symmetry.
\newblock {\em Phys. Lett. A}, 416:Paper No. 127669, 6, 2021.

\bibitem[FAL{\etalchar{+}}24]{FarAkhLawBiaSir_pre24}
O.~Farooq, A.~Akhshani, M.~\L{}awniczak, M.~Bia\l{}ous, and L.~Sirko.
\newblock Coupled unidirectional chaotic microwave graphs.
\newblock {\em Phys. Rev. E}, 110:014206, Jul 2024.

\bibitem[GA04]{GnuAlt_prl04}
S~Gnutzmann and A~Altland.
\newblock Universal spectral statistics in quantum graphs.
\newblock {\em Phys. Rev. Lett.}, 93:194101, 2004.

\bibitem[GA05]{GnuAlt_pre05}
S.~Gnutzmann and A.~Altland.
\newblock Spectral correlations of individual quantum graphs.
\newblock {\em Phys. Rev. E (3)}, 72(5):056215, 14, 2005.

\bibitem[GO13a]{GutOsi_jsp13}
B.~Gutkin and V.~Osipov.
\newblock Clustering of periodic orbits and ensembles of truncated unitary
  matrices.
\newblock {\em J. Stat. Phys.}, 153(6):1049--1064, 2013.

\bibitem[GO13b]{GutOsi_nonlin13}
B.~Gutkin and V.~Osipov.
\newblock Clustering of periodic orbits in chaotic systems.
\newblock {\em Nonlinearity}, 26(1):177--200, 2013.

\bibitem[GS06]{GnuSmi_ap06}
S.~Gnutzmann and U.~Smilansky.
\newblock Quantum graphs: Applications to quantum chaos and universal spectral
  statistics.
\newblock {\em Adv. Phys.}, 55(5--6):527--625, 2006.

\bibitem[GS07]{GavSmi_jpa07}
U.~Gavish and U.~Smilansky.
\newblock Degeneracies in the length spectra of metric graphs.
\newblock {\em J. Phys. A}, 40(33):10009--10020, 2007.

\bibitem[GS24]{GnuSmi_jpa24}
S.~Gnutzmann and U.~Smilansky.
\newblock Information scrambling and chaos induced by a {H}ermitian matrix.
\newblock {\em J. Phys. A}, 57(37):Paper No. 37LT01, 13, 2024.

\bibitem[HBPan{\etalchar{+}}04]{HulBauPakSavZycSir_pre04}
O.~Hul, S.~Bauch, P.~Pako\ifmmode~\acute{n}\else \'{n}\fi{}ski, N.~Savytskyy,
  K.~\ifmmode~\dot{Z}\else \.{Z}\fi{}yczkowski, and L.~Sirko.
\newblock Experimental simulation of quantum graphs by microwave networks.
\newblock {\em Phys. Rev. E}, 69:056205, May 2004.

\bibitem[HH22a]{HarHud_jpa22}
J.~M. Harrison and T.~Hudgins.
\newblock Complete dynamical evaluation of the characteristic polynomial of
  binary quantum graphs.
\newblock {\em J. Phys. A}, 55(42):Paper No. 425202, 45, 2022.

\bibitem[HH22b]{HarHud_epl22}
J.~M. Harrison and T.~Hudgins.
\newblock Periodic orbit evaluation of a spectral statistic of quantum graphs
  without the semiclassical limit.
\newblock {\em Europhysics Letters}, 138(3):30002, may 2022.

\bibitem[HLKS21]{HoLuKuSto_pre21}
T.~Hofmann, J.~Lu, U.~Kuhl, and H.-J. St\"ockmann.
\newblock Spectral duality in graphs and microwave networks.
\newblock {\em Phys. Rev. E}, 104:045211, Oct 2021.

\bibitem[HS11]{OleSir_pre11}
O.~Hul and L.~Sirko.
\newblock Parameter-dependent spectral statistics of chaotic quantum graphs:
  Neumann versus circular orthogonal ensemble boundary conditions.
\newblock {\em Phys. Rev. E}, 83:066204, Jun 2011.

\bibitem[HS19]{HarSwi_jpa19}
J.~M. Harrison and E.~Swindle.
\newblock Spectral properties of quantum circulant graphs.
\newblock {\em J. Phys. A}, 52(28):285101, 33, 2019.

\bibitem[HW12]{HarWin_jpa12}
J.~M. Harrison and B.~Winn.
\newblock Intermediate statistics for a system with symplectic symmetry: the
  {D}irac rose graph.
\newblock {\em J. Phys. A}, 45(43):435101, 23, 2012.

\bibitem[JMS14]{JoyMueSie_epl14}
C.~H. Joyner, S.~Mueller, and M.~Sieber.
\newblock {GSE} statistics without spin.
\newblock {\em EPL}, 107(5):50004, 2014.

\bibitem[KN23]{KosNic_book22}
A.~Kostenko and N.~Nicolussi.
\newblock {\em Laplacians on infinite graphs}, volume~3 of {\em Memoirs of the
  European Mathematical Society}.
\newblock EMS Press, Berlin, 2022 \copyright 2023.

\bibitem[KS97]{KotSmi_prl97}
T.~Kottos and U.~Smilansky.
\newblock Quantum chaos on graphs.
\newblock {\em Phys. Rev. Lett.}, 79(24):4794--4797, 1997.

\bibitem[KS99a]{KosSch_jpa99}
V.~Kostrykin and R.~Schrader.
\newblock Kirchhoff's rule for quantum wires.
\newblock {\em J. Phys. A}, 32(4):595--630, 1999.

\bibitem[KS99b]{KotSmi_ap99}
T.~Kottos and U.~Smilansky.
\newblock Periodic orbit theory and spectral statistics for quantum graphs.
\newblock {\em Ann. Physics}, 274(1):76--124, 1999.

\bibitem[KS01]{KotSch_phys01}
T.~Kottos and H.~Schanz.
\newblock Quantum graphs: a model for quantum chaos.
\newblock {\em Physica E: Low-dimensional Systems and Nanostructures},
  9(3):523--530, 2001.
\newblock Proceedings of an International Workshop and Seminar on the Dynamics
  of Complex Systems.

\bibitem[Kur24]{Kurasov_book24}
P.~Kurasov.
\newblock {\em Spectral geometry of graphs}, volume 293 of {\em Operator
  Theory: Advances and Applications}.
\newblock Birkh\"auser/Springer, Berlin, [2024] \copyright 2024.

\bibitem[LHKS23]{LuHofKuhSto_ent23}
J.~Lu, T.~Hofmann, U.~Kuhl, and H.-J. St\"{o}ckmann.
\newblock Implications of spectral interlacing for quantum graphs.
\newblock {\em Entropy}, 25(1), 2023.

\bibitem[Luo25]{Luo_arxiv25}
Ye~Luo.
\newblock On a trace formula of counting eulerian cycles.
\newblock {\em arXiv:2502.02915}, 2025.

\bibitem[MR98]{McKayRob_cpc98}
B.~D. McKay and R.~W. Robinson.
\newblock Asymptotic enumeration of eulerian circuits in the complete graph.
\newblock {\em Combinatorics, Probability and Computing}, 7(4):437–449, 1998.

\bibitem[PcvW13]{PluWei_prl13}
Z.~Pluha\ifmmode~\check{r}\else \v{r}\fi{} and H.~A. Weidenm\"uller.
\newblock Universal chaotic scattering on quantum graphs.
\newblock {\em Phys. Rev. Lett.}, 110:034101, Jan 2013.

\bibitem[PcvW14]{PluWei_prl14}
Z.~Pluha\ifmmode~\check{r}\else \v{r}\fi{} and H.~A. Weidenm\"uller.
\newblock Universal quantum graphs.
\newblock {\em Phys. Rev. Lett.}, 112:144102, Apr 2014.

\bibitem[PcvW15]{PluWei_jpa15}
Z.~Pluha\ifmmode~\check{r}\else \v{r}\fi{} and H.~A. Weidenm\"uller.
\newblock Quantum graphs and random-matrix theory.
\newblock {\em J. Phys. A}, 48(27):275102, 30, 2015.

\bibitem[RAJ{\etalchar{+}}16]{RehAllJoyMueSieKuhSto_prl16}
A.~Rehemanjiang, M.~Allgaier, C.~H. Joyner, S.~M\"uller, M.~Sieber, U.~Kuhl,
  and H.-J. St\"ockmann.
\newblock Microwave realization of the gaussian symplectic ensemble.
\newblock {\em Phys. Rev. Lett.}, 117:064101, Aug 2016.

\bibitem[RB69]{RofeBeketov1969}
F.~S. Rofe-Beketov.
\newblock Selfadjoint extensions of differential operators in a space of
  vector-valued functions.
\newblock {\em Teor. Funkci\u i\ Funkcional. Anal. i Prilo\v zen.}, (8):3--24,
  1969.

\bibitem[SS00]{SchSmi_pmb00}
H.~Schanz and U.~Smilansky.
\newblock Spectral statistics for quantum graphs: periodic orbits and
  combinatorics.
\newblock {\em Phil. Mag. B}, 80(12):1999--2021, 2000.

\bibitem[SS01]{SchSmi_ejc01}
H.~Schanz and U.~Smilansky.
\newblock Combinatorial identities from the spectral theory of quantum graphs.
\newblock {\em Electron. J. Combin.}, 8(2):Research Paper 16, 16 pp.
  (electronic), 2001.
\newblock In honor of Aviezri Fraenkel on the occasion of his 70th birthday.

\bibitem[ST41]{SmiTut_amm1941}
C.~A.~B. Smith and W.~T. Tutte.
\newblock On unicursal paths in a network of degree 4.
\newblock {\em American Mathematical Monthly}, 48:233--237, 1941.

\bibitem[Tan00]{Tan_jpa00}
G.~Tanner.
\newblock Spectral statistics for unitary transfer matrices of binary graphs.
\newblock {\em J. Phys. A}, 33(18):3567--3585, 2000.

\bibitem[Tan01]{Tan_jpa01}
G.~Tanner.
\newblock Unitary-stochastic matrix ensembles and spectral statistics.
\newblock {\em J. Phys. A}, 34(41):8485--8500, 2001.

\bibitem[vAEdB51]{EhrBru_sis1951}
T.~van Aardenne-Ehrenfest and N.~G. de~Bruijn.
\newblock Circuits and trees in oriented linear graphs.
\newblock {\em Simon Stevin}, 28:203--217, 1951.

\bibitem[Wei20]{Weidenmueller_chapter20}
H.~A. Weidenm{\"u}ller.
\newblock {\em Massive Modes for Quantum Graphs}, pages 341--357.
\newblock Springer International Publishing, Cham, 2020.

\end{thebibliography}

\end{document}